\documentclass{article}
\usepackage{makeidx}
\usepackage{amssymb}
\usepackage{amsfonts}
\usepackage{amsmath}
\usepackage[backref,colorlinks=true]{hyperref}

\newtheorem{theorem}{Theorem}

\newtheorem{corollary}[theorem]{Corollary}

\newtheorem{lemma}[theorem]{Lemma}

\begin{document}

\title{Properties of multitype subcritical branching processes in random
environment \thanks{%
This work was supported by the Russian Science Foundation under the grant
17-11-01173-Ext}}
\author{Vladimir Vatutin\thanks{%
Novosibirsk State University, Novosibirsk, Russia, e-mail: vatutin@mi-ras.ru}%
, Elena Dyakonova\thanks{%
Novosibirsk State University, Novosibirsk, Russia, e-mail: elena@mi-ras.ru}}

\date{}
\maketitle

\begin{abstract}
We study properties of a $p-$type subcritical branching process in random
environment initiated at moment zero by a vector $\mathbf{z}=\left(
z_{1},..,z_{p}\right) $\ of particles of different types. Assuming that the
process belongs to the class of the so-called strongly subcritical processes
we show that its survival probability to moment $n$\ behaves for large $n$\
as $C(\mathbf{z})\lambda ^{n}$\ where $\lambda $\ is the upper Lyapunov
exponent for the product of mean matrices of the process and $C(\mathbf{z})$%
\ is an explicitly given constant. We also demonstrate that the limiting
conditional distribution of the number of particles given the survival of
the process for a long time does not depend on the vector $\mathbf{z}$ of
the number of particles initiated the process.
\end{abstract}



\section{Introduction and statement of main results}

We investigate asymptotic properties of a class of subcritical multitype
branching processes in random environment (MBPRE's). To formulate the
results of the paper we introduce some notation.

For $p$-dimensional vectors $\mathbf{x}=\left( x_{1},\ldots ,x_{p}\right) $
and $\mathbf{y}=\left( y_{1},\ldots ,y_{p}\right) $ we set
\begin{equation*}
(\mathbf{x},\mathbf{y}):=\sum_{i=1}^{p}x_{i}y_{i},\quad |\mathbf{x}%
|:=\sum_{i=1}^{p}|x_{i}|\quad \text{and}\quad \mathbf{x}^{\mathbf{y}%
}:=\prod_{i=1}^{p}(x_{i})^{y_{i}}.
\end{equation*}%
The standard basis vectors will be denoted by $\mathbf{e}_{i}$, $%
i=1,2,\ldots ,p.$ We set $\mathbf{0}:=(0,0,\ldots ,0)$ and $\mathbf{1}%
:=(1,1,\ldots ,1)$, denote by $\mathbb{N}_{0}^{p}$ the set of all $p$%
-dimensional vectors with non-negative integer-valued components and
introduce the notation $\mathbb{N}_{+}^{p}:=\mathbb{N}_{0}^{p}\backslash
\left\{ \mathbf{0}\right\} $.

For every $p$-tuple $(\mu ^{1},\mu ^{2},\ldots ,\mu ^{p})$ of probability
measures on $\mathbb{N}_{0}^{p},$ we define its multidimensional
multivariate generating function $\mathbf{f}=(f^{1},f^{2},\ldots ,f^{p})$ by
the relations
\begin{equation*}
f^{i}(\mathbf{s}):=\sum_{\mathbf{z}\in \mathbb{N}_{0}^{p}}\mu ^{i}\{\mathbf{z%
}\}\mathbf{s}^{\mathbf{z}},\quad \mathbf{s=}\left( s_{1},...,s_{p}\right)
\in \lbrack 0,1]^{p},\ i=1,2,\ldots ,p.
\end{equation*}%
Any sequence $\{\mathbf{f}_{n}=(f_{n}^{1},f_{n}^{2},\ldots ,f_{n}^{p}),n\geq
1\}$ of multidimensional multivariate generating functions is called a
\textit{varying environment}. The corresponding measure $\mu _{n}^{i}$
describes the offspring law of type $i$ particles in generation $n-1$.

The model we are planning to investigate is a $p$-type branching process $%
\mathbf{Z}_{n}=(Z_{n}^{1},Z_{n}^{2},\ldots ,Z_{n}^{p})$, $n\geq 0,$ where $%
Z_{n}^{i}$ is the number of type $i$ particles in the process at moment $n$.
This process has a (may be random) starting point $\mathbf{Z}_{0}$ and the
population sizes of the subsequent generations of the process are specified
by the recursion
\begin{equation}
Z_{n}^{i}=\sum_{j=1}^{p}\sum_{k=1}^{Z_{n-1}^{j}}X_{n,j,k}^{i},\quad
i=1,2,\ldots ,p,\ n\geq 1,  \label{def.varying}
\end{equation}%
where $\mathbf{X}_{n,j,k}=(X_{n,j,k}^{1},X_{n,j,k}^{2},\ldots
,X_{n,j,k}^{p}) $, ${k\geq 1,}$ are independent random vectors distributed
according to $\mu _{n}^{j}$. Here and in what follows we denote random
objects by upper case symbols using the respective bold symbols if the
objects are vectors or matrices.

Equipping the set of all tuples $(\mu ^{1},\mu ^{2},\ldots ,\mu ^{p})$ with
the total variation distance we obtain a metric space $\mathcal{P}_{p}(%
\mathbb{N}_{0}^{p})$ which, due to the one-to-one correspondence $(\mu
^{1},\mu ^{2},\ldots ,\mu ^{p})\leftrightarrows \mathbf{f}%
=(f^{1},f^{2},\ldots ,f^{p})$, we identify with the space of all
multidimensional multivariate probability generating functions. Thus, we are
able to specify a probability measure $\mathbb{P}$ on the space of such
generating functions. This agreement allows us to introduce a sequence of
independent, identically distributed random $p-$dimensional generating
functions $\mathcal{E}:=\{\mathbf{F}_{n}=(F_{n}^{1},F_{n}^{2},\ldots
,F_{n}^{p}),n\geq 1\}$ which will be called a \textit{random environment}.
We say that $\mathcal{Z}=\{\mathbf{Z}_{n},n\geq 0\}$ is a $p-$type BPRE, if
for every fixed realisation of the environmental sequence its conditional
distribution is determined by \eqref{def.varying}.

BPRE's with one type of particles have been intensively investigated during
the last two decades and many their properties are well understood. The
reader may find a unified presentation of the corresponding results in \cite%
{KV2017}. The multi-dimensional case is much less studied. For instance,
only recently estimates for the survival probability of the critical and
some classes of subcritical multitype BPRE's were found under relatively
general conditions, see \cite{Dyak11},\cite{Dyak12},\cite{Dy2013},\cite%
{Dyak15},\cite{LPP2016}, \cite{VD2017}, \cite{VW2018}, \cite{VD2019}.

It known that asymptotic properties of MBPRE's are often described in terms
of properties of the (random) mean matrices
\begin{equation*}
\mathbf{M}_{n}=\left( M_{ij}(n)\right) _{i,j=1}^{p}:=\left( \frac{\partial
F_{n}^{i}}{\partial s_{j}}(\mathbf{1})\right) _{i,j=1}^{p},\quad n\geq 1.
\end{equation*}%
If $\mathbf{F}_{n}$ are independent and distributed as a generating function
$\mathbf{F}=\left( F^{1},\ldots ,F^{p}\right) $ then $\mathbf{M}_{n}$ are
independent probabilistic copies of the random matrix
\begin{equation*}
\mathbf{M=}\left( M_{ij}\right) _{i,j=1}^{p}:=\left( \frac{\partial F^{i}}{%
\partial s_{j}}(\mathbf{1})\right) _{i,j=1}^{p}.
\end{equation*}%
We assume that the distribution of $\mathbf{M}$ satisfies the following
assumptions:

\textbf{Condition }$\mathbf{H1}$. The set $\Theta _{+}:=\left\{ \theta >0:%
\mathbb{E}\left[ \left\Vert \mathbf{M}\right\Vert ^{\theta }\right] <\infty
\right\} $ is nonempty, where $\Vert \mathbf{M}\Vert $ is the operator norm
of $\mathbf{M}$.

\textbf{Condition }$\mathbf{H2}$. The support of the distribution of $%
\mathbf{M}$ acts strongly irreducibly on the semi-group of matrices with
non-negative entries, i.e. no proper finite union of subspaces of $\mathbb{R}%
^{p}$ is invariant with respect to all elements of the multiplicative
semi-group $\mathbb{S}^{+}$ of $p\times p$ matrices generated by the support
of $\mathbf{M}$.

\textbf{Condition }$\mathbf{H3}$. There exists a positive number $\gamma >1$
such that
\begin{equation*}
1\leq \frac{\max_{i,j}M_{ij}}{\min_{i,j}M_{ij}}\leq \gamma .
\end{equation*}

To introduce one more condition we define the cone
\begin{equation*}
\mathcal{C=}\left\{ \mathbf{x}=(x_{1},...,x_{p})\in \mathbb{R}^{p}:x_{i}\geq
0\text{ for any }i=1,...,p\right\} ,
\end{equation*}%
and the space
\begin{equation*}
\mathbb{X}=\mathcal{C}\cap \left\{ \mathbf{x}:\mathbf{x}\in \mathbb{R}%
^{p},\left\vert \mathbf{x}\right\vert =1\right\}.
\end{equation*}

\textbf{Condition }$\mathbf{H4}$. There exists $\delta >0$ such that

\begin{equation*}
\mathbb{P}\left( \mathbf{h}\in \mathbb{S}^{+}:\text{ for any }\mathbf{x}\in
\mathbb{X},\,\log \left\vert \mathbf{xh}\right\vert \geq \delta \right) >0.
\end{equation*}

To formulate the next assumption we define matrices
\begin{equation*}
\mathbf{B}(k):=\left( \frac{\partial ^{2}F^{k}}{\partial s_{i}\partial s_{j}}%
(\mathbf{1})\right) _{i,j=1}^{p}\text{ and }\mathbf{B}_{n}(k):=\left( \frac{%
\partial ^{2}F_{n}^{k}}{\partial s_{i}\partial s_{j}}(\mathbf{1})\right)
_{i,j=1}^{p}
\end{equation*}%
and random variables
\begin{equation*}
\mathcal{T}:=\frac{1}{\Vert \mathbf{M}\Vert ^{2}}\sum_{k=1}^{p}\Vert \mathbf{%
B}(k)\Vert \text{ and }\mathcal{T}_{n}:=\frac{1}{\Vert \mathbf{M}_{n}\Vert
^{2}}\sum_{k=1}^{p}\Vert \mathbf{B}_{n}(k)\Vert , n=1,2,\ldots .
\end{equation*}%
Thus, $\mathcal{T}_{n}$ are independent probabilistic copies of $\mathcal{T}%
. $

\textbf{Condition }$\mathbf{H5}$. There exists $\varepsilon >0$ such that
\begin{equation*}
\mathbb{E}\left[ \Vert \mathbf{M}\Vert |\log \mathcal{T}|^{1+\varepsilon }%
\right] <\infty .
\end{equation*}

Conditions $\mathbf{H1}-\mathbf{H5}$\ impose restrictions on the integral
characteristics of the MPBRE. Our next assumption directly concerns the
reproduction laws of particles and forces particles of all types to die
without children with positive probability.

\textbf{Condition }$\mathbf{H6}$\textbf{. }$\mathbb{P}(F^{k}[\mathbf{0}%
]>0,k=1,...,p)=1$.

It is well known (see, for instance,\ point 2.1.2 in \cite{BDGM2014}) that
the limit
\begin{equation*}
\lambda \left( \theta \right) :=\lim_{n\rightarrow \infty }\left( \mathbb{E}%
\left[ \left\Vert \mathbf{M}_{n}\cdot \cdot \cdot \mathbf{M}_{1}\right\Vert
^{\theta }\right] \right) ^{1/n}<\infty
\end{equation*}%
exists for every $\theta \in \Theta :=\left\{ 0\right\} \cup \Theta _{+}$.
Set
\begin{equation*}
\Lambda (\theta ):=\log \lambda (\theta ),\quad \theta \in \Theta .
\end{equation*}%
Along with the random mean matrix $\mathbf{M}$ we consider its expectation
\begin{equation*}
\mathbf{m}:=\mathbb{E}\mathbf{M}=\left( \mathbb{E}M_{ij}\right)
_{i,j=1}^{p}=\left( m_{ij}\right) _{i,j=1}^{p}.
\end{equation*}%
Clearly, the elements of $\mathbf{m}$ are positive. Let $\lambda $ be the
Perron root of $\mathbf{m}$ and let $\mathbf{V}=\left( V_{1},..,V_{p}\right)
$ and $\mathbf{U}=\left( U_{1},..,U_{p}\right) $ be the strictly positive
left and right eigenvectors of $\mathbf{m}$ corresponding to the eigen-value
$\lambda $
\begin{equation}
\mathbf{Vm=}\lambda \mathbf{V,\ mU=}\lambda \mathbf{U}\text{\ and\textbf{\ }%
scaled as\textbf{\ }}(\mathbf{V,1})=1,\ (\mathbf{V,U})=1\text{.}
\label{MuEigen}
\end{equation}

In view of (\ref{MuEigen}) and positivity of the components of $\mathbf{U}$%
\begin{equation*}
\lim_{n\rightarrow \infty }\left\vert \mathbf{m}^{n}\right\vert
^{1/n}=\lim_{n\rightarrow \infty }\left\vert \mathbf{m}^{n}\mathbf{U}%
\right\vert ^{1/n}=\lambda .
\end{equation*}%
Observe further that
\begin{eqnarray*}
\frac{1}{p\left\vert \mathbf{U}\right\vert }\left\vert \mathbf{M}_{n}\cdot
\cdot \cdot \mathbf{M}_{1}\mathbf{U}\right\vert &\leq &\left\Vert \mathbf{M}%
_{n}\cdot \cdot \cdot \mathbf{M}_{1}\right\Vert =\sup_{\left\vert \mathbf{x}%
\right\vert =1}\left\vert \mathbf{M}_{n}\cdot \cdot \cdot \mathbf{M}_{1}%
\mathbf{x}\right\vert \\
&\leq &\left\vert \mathbf{M}_{n}\cdot \cdot \cdot \mathbf{M}_{1}\mathbf{1}%
\right\vert \leq \frac{1}{\min_{1\leq i\leq p}U_{i}}\left\vert \mathbf{M}%
_{n}\cdot \cdot \cdot \mathbf{M}_{1}\mathbf{U}\right\vert .
\end{eqnarray*}%
Thus, by independency of random matrices and the Perron-Frobenius theorem
\begin{eqnarray*}
\lambda \left( 1\right) &=&\lim_{n\rightarrow \infty }\left( \mathbb{E}\left[
\left\Vert \mathbf{M}_{n}\cdot \cdot \cdot \mathbf{M}_{1}\right\Vert \right]
\right) ^{1/n} \\
&=&\lim_{n\rightarrow \infty }\left( |\left( \mathbb{E}\mathbf{M}_{n}\cdot
\cdot \cdot \mathbb{E}\mathbf{M}_{1}\right) \mathbf{U}|\right)
^{1/n}=\lim_{n\rightarrow \infty }\left\vert \mathbf{m}^{n}\mathbf{U}%
\right\vert ^{1/n}=\lambda .
\end{eqnarray*}%
Therefore, $\lambda \left( 1\right) $ is the Perron~root of $\mathbf{m}$.
For this reason we rather often write below $\lambda $ for $\lambda \left(
1\right) $.

Now we are ready to formulate the main results of the paper describing
properties of a MBPRE which is initiated by a deterministic vector
\begin{equation*}
\mathbf{Z}_{0}=(Z_{0}^{1},Z_{0}^{2},\ldots ,Z_{0}^{p})=\mathbf{z}=\left(
z_{1},z_{2,}...,z_{p}\right) \in \mathbb{N}_{+}^{p}
\end{equation*}%
of particles at time $0$. Introduce the event
\begin{equation}
\mathcal{D}_{n}\left( \mathbf{z}\right) :=\left\{ \mathbf{Z}_{n}\neq \mathbf{%
0}\big|\mathbf{Z}_{0}=\mathbf{z}\right\}  \label{DefB}
\end{equation}%
and denote by $\mathrm{int}\Theta $ the interior of the set $\Theta $.

\begin{theorem}
\label{T_condSeveral} Let Conditions $\mathbf{H1}-\mathbf{H5}$ be valid, $%
1\in \mathrm{int}\Theta $ and \ $\Lambda ^{\prime }(1)<0$. Then for all%
\textbf{\ }$\mathbf{z}\in N_{+}^{p}$,\ as $n\rightarrow \infty ,$%
\begin{equation}
\mathbb{P}\left( \mathcal{D}_{n}\left( \mathbf{z}\right) \right) \sim \frac{(%
\mathbf{z},\mathbf{U})}{\mathbb{E}\left[ \left( \mathbf{T}_{\infty },\mathbf{%
U}\right) \right] }\lambda ^{n}(1)  \label{SinglStrong}
\end{equation}%
and, for any \ $\mathbf{s\in }\left[ 0,1\right] ^{p}$ and all $\mathbf{z\in }%
\mathbb{N}_{+}^{p}$
\begin{equation}
\lim_{n\rightarrow \infty }\mathbb{E}\left[ \mathbf{s}^{\mathbf{Z}_{n}}|%
\mathbf{Z}_{n}\neq \mathbf{0};\mathbf{Z}_{0}=\mathbf{z}\right] =\sum_{%
\mathbf{x}\in \mathbb{N}_{+}^{p}}t_{\mathbf{x}}\mathbf{s}^{\mathbf{x}%
}:=T\left( \mathbf{s}\right) =\mathbb{E}\left[ \mathbf{s}^{\mathbf{T}%
_{\infty }}\right]  \label{EqYz}
\end{equation}%
exists, is independent of $\mathbf{z}$, specifies a proper distribution on $%
\mathbb{N}_{+}^{p}$ and solves the equation
\begin{equation}
\mathbb{E}\left[ T\left( \mathbf{F}(\mathbf{s})\right) \right] =\lambda
(1)T\left( \mathbf{s}\right) +1-\lambda (1).  \label{EqY}
\end{equation}
\end{theorem}

\textbf{Remark 1}. The case $|\mathbf{z}|=1$ was considered in \cite{VW2018}%
, where, for all $i=1,...,p$ it was shown that%
\begin{equation}
\mathbb{P}\left( \mathcal{D}_{n}\left( \mathbf{e}_{i}\right) \right) \sim K(%
\mathbf{e}_{i})\lambda ^{n}(1)  \label{SinglStrong0}
\end{equation}%
as $n\rightarrow \infty $ and that
\begin{equation}
\lim_{n\rightarrow \infty }\mathbb{E}\left[ \mathbf{s}^{\mathbf{Z}_{n}}|%
\mathbf{Z}_{n}\neq \mathbf{0};\mathbf{Z}_{0}=\mathbf{e}_{i}\right] =T_{i}(%
\mathbf{s})  \label{EqYz0}
\end{equation}%
exists. However, no explicit expressions for the constants $K(\mathbf{e}%
_{i}) $ in (\ref{SinglStrong0}) were given and independence of the limit in (%
\ref{EqYz0}) on the initial value $\mathbf{z}=\mathbf{e}_{i},i=1,...,p$ was
not established. Our Theorem \ref{T_condSeveral} fills this gap. The proofs
of (\ref{SinglStrong}) and (\ref{EqYz}) are based on a combination of the
methods from \cite{VW2018} used in Section \ref{Sec1} and the change of
measure of a new type applied in Section \ref{SecPrecise}.

\textbf{Remark 2}. If $\Lambda ^{\prime }(1)<0$ then $\Lambda ^{\prime
}(0)<0 $ by convexity of $\lambda (\theta )$. In particular, if $p=1$ then
the inequalities $\Lambda ^{\prime }(0)<0$ and $\Lambda ^{\prime }(1)<0$ are
reduced to $\mathbb{E}[M_{11}]<0$ and $\mathbb{E}[M_{11}\log M_{11}]<0$.
Such single-type BPRE's (called strongly subcritical) were investigated, for
instance, in \cite{agkv2}, \cite{abkv}, \cite{ABKV2014}, \cite{GKV2003} and
\cite{GL2001}.

Multitype subcritical branching processes were considered in \cite{Dyak07}
and \cite{Dy2013} under the additional restriction that all mean matrices of
the reproduction laws have a common deterministic left or right eigenvector
corresponding to the Perron roots of these matrices. We do not require the
validity of such condition.

The following two corollaries provide an additional information about the
limiting distributions of the process $\left\{ \mathbf{Z}_{n},n\geq
1\right\} $ under two types of conditioning.

\begin{corollary}
\label{C_Qprocess2}Let the conditions of Theorem \ref{T_condSeveral} be
valid. Then

(i) for each fixed $m\geq 1$ and $\mathbf{z,j}_{0},\mathbf{j}_{1},...,%
\mathbf{j}_{m}\in \mathbb{N}_{+}^{p}$
\begin{eqnarray}
&&\lim_{n\rightarrow \infty }\mathbb{P}\left( \mathbf{Z}_{n-m}=\mathbf{j}%
_{0},\mathbf{Z}_{n-m+1}=\mathbf{j}_{1},...,\mathbf{Z}_{n}=\mathbf{j}_{m}|%
\mathbf{Z}_{n}\neq \mathbf{0;Z}_{0}=\mathbf{z}\right)  \notag \\
&=&\frac{1}{\lambda ^{m}(1)}\mathbb{P}\left( \mathbf{Z}_{1}=\mathbf{j}%
_{1},...,\mathbf{Z}_{m}=\mathbf{j}_{m}|\mathbf{Z}_{m}\neq \mathbf{0;Z}_{0}=%
\mathbf{j}_{0}\right) \mathbb{P}\left( \mathbf{T}_{\infty }=\mathbf{j}%
_{m}\right) ;  \label{PrelMult}
\end{eqnarray}%
(ii) if the tuple $0=n_{0}<n_{1}<....<n_{r}=n$ is such that
\begin{equation*}
n_{\ast }:=\min_{0\leq k\leq r-1}\left( n_{k+1}-n_{k}\right) \rightarrow
\infty ,
\end{equation*}%
then%
\begin{eqnarray}
&&\mathbb{P}(\mathbf{Z}_{n_{1}}=\mathbf{j}_{1},...,\mathbf{Z}_{n_{r}}=%
\mathbf{j}_{r}|\mathbf{Z}_{n}\neq \mathbf{0;Z}_{0}=\mathbf{z})  \notag \\
&&\qquad \qquad \qquad \rightarrow \left( \prod_{i=1}^{r-1}\left( \mathbf{j}%
_{i},\mathbf{K}\right) \mathbb{P}\left( \mathbf{T}_{\infty }=\mathbf{j}%
_{i}\right) \right) \mathbb{P}\left( \mathbf{T}_{\infty }=\mathbf{j}%
_{r}\right) ,  \label{LimitMult}
\end{eqnarray}%
where%
\begin{equation}
\mathbf{K}=\frac{\mathbf{U}}{\mathbb{E}\left[ \left( \mathbf{T}_{\infty },%
\mathbf{U}\right) \right] }.  \label{DefK}
\end{equation}
\end{corollary}

Thus, given survival of the strongly subcritical MBPRE to a distant moment $%
n $, the vectors of the number of particles at moments being far from each
other are asymptotically independent. Such phenomenon for single-type
strongly subcritical BPRE's was described in \cite{agkv2}.

We now introduce the so-called $Q$-process $\mathcal{Y}:=\left\{ \mathbf{%
\hat{Y}}_{n},n\geq 1\right\} $ which may be nonrigorously considered as the
process $\left\{ \mathbf{Z}_{n},n\geq 0\right\} $ conditioned to survive in
the distant future and whose multi-dimensional distributions are specified
by the formulas%
\begin{eqnarray*}
&&\mathbb{P}\left( \mathbf{\hat{Y}}_{1}=\mathbf{j}_{1},\mathbf{\hat{Y}}_{2}=%
\mathbf{j}_{2},...,\mathbf{\hat{Y}}_{n}=\mathbf{j}_{n}|\mathbf{\hat{Y}}_{0}=%
\mathbf{y}\right) \\
&&\qquad =\lim_{m\rightarrow \infty }\mathbb{P}\left( \mathbf{Z}_{1}=\mathbf{%
j}_{1},\mathbf{Z}_{2}=\mathbf{j}_{2},...,\mathbf{Z}_{n}=\mathbf{j}_{n}|%
\mathbf{Z}_{n+m}\neq \mathbf{0;Z}_{0}=\mathbf{y}\right) .
\end{eqnarray*}

The next statement gives a more explicit representation for the
distributions of $\mathcal{Y}$.

\begin{corollary}
\label{C_Qprocess}Under the conditions of Theorem \ref{T_condSeveral} for
each fixed $n\geq 1$
\begin{eqnarray*}
&&\mathbb{P}\left( \mathbf{\hat{Y}}_{1}=\mathbf{j}_{1},\mathbf{\hat{Y}}_{2}=%
\mathbf{j}_{2},...,\mathbf{\hat{Y}}_{n}=\mathbf{j}_{n}|\mathbf{\hat{Y}}_{0}=%
\mathbf{y}\right) \\
&&\qquad=\frac{1}{\lambda ^{n}(1)}\frac{\left( \mathbf{j}_{n},\mathbf{K}%
\right) }{(\mathbf{y},\mathbf{K})}\mathbb{P}\left( \mathbf{Z}_{1}=\mathbf{j}%
_{1},\mathbf{Z}_{2}=\mathbf{j}_{2},...,\mathbf{Z}_{n}=\mathbf{j}_{n}\mathbf{%
|Z}_{0}=\mathbf{y}\right)
\end{eqnarray*}%
and for any $\mathbf{j},\mathbf{y}\in \mathbb{N}_{+}^{p}$
\begin{equation*}
\lim_{n\rightarrow \infty }\mathbb{P}\left( \mathbf{\hat{Y}}_{n}=\mathbf{j}|%
\mathbf{\hat{Y}}_{0}=\mathbf{y}\right) =\frac{\left( \mathbf{j,U}\right) }{%
\mathbb{E}\left[ \left( \mathbf{T}_{\infty },\mathbf{U}\right) \right] }%
\mathbb{P}\left( \mathbf{T}_{\infty }=\mathbf{j}\right) .
\end{equation*}
\end{corollary}

Thus, the $Q$-process $\mathcal{Y}$ has a kind of size-biased distribution
in the limit.

The rest of the paper is organized as follows. To prove the desired
statements we need to perform two changes of measure. The first change of
measure is introduced in Section \ref{Sec1}. Using the new measure we obtain
in Section \ref{Sec2} a rough representation for the survival probability of
a class of subcritical MBPRE's and prove in Section \ref{Sec3} a Yaglom-type
limit theorem for the distribution of the number of particles of different
types given the process is initiated by several particles. In Section \ref%
{SecPrecise} we make the second change of measure and complete the proof of
Theorem \ref{T_condSeveral} by justifying asymptotic representation (\ref%
{SinglStrong}). A description of properties of $Q$-processes is given in
Section~\ref{SecQproc}.

\section{First change of measure\label{Sec1}}

We agree to denote $p\times p$ deterministic matrices with non-negative
entries by bold lower case symbols and, as a rule, use\textit{\ the same
notation for row and column vectors}. For instance, if $\mathbf{x}$ is a $p$%
-dimensional vector and $\mathbf{h}$ is a $p\times p$ deterministic matrix
then, the $\mathbf{x}$ will be treated as a row vector in the product $%
\mathbf{xh}$ and as a column vector in the product $\mathbf{hx}$. It will be
clear from the context which form is used.

For every $\mathbf{x\in }\mathbb{X}$ and a $p\times p$ matrix $\mathbf{h\in }%
\mathbb{S}^{+}$ we specify the column vector
\begin{equation}
\mathbf{h}\circ \mathbf{x}:=\frac{\mathbf{hx}}{|\mathbf{hx}|}\in \mathbb{X}%
\text{ }  \label{Column}
\end{equation}%
if $\mathbf{hx}\neq \mathbf{0}$ and the row vector%
\begin{equation}
\text{ }\mathbf{x\circ h}:=\frac{\mathbf{xh}}{|\mathbf{xh}|}\in \mathbb{X}%
\text{ }  \label{Row}
\end{equation}%
if $\mathbf{xh}\neq \mathbf{0}$.

Denote by $\mathcal{C}_{b}\left( \mathbb{X}\right) $ the set of all bounded
continuous functions on $\mathbb{X}$. For $\theta \in \Theta ,$ $g\in
\mathcal{C}_{b}\left( \mathbb{X}\right) ,$ and $\mathbf{x\in }\mathbb{X}$
define the transition operators
\begin{equation*}
P_{\theta }g(\mathbf{x}):=\mathbb{E}\left[ \left\vert \mathbf{Mx}\right\vert
^{\theta }g\left( \mathbf{M\circ x}\right) \right]
\end{equation*}%
and
\begin{equation*}
P_{\theta }^{\ast }g(\mathbf{x}):=\mathbb{E}\left[ \left\vert \mathbf{M}^{T}%
\mathbf{x}\right\vert ^{\theta }g\left( \mathbf{M}^{T}\mathbf{\circ x}%
\right) \right] .
\end{equation*}

If the assumptions $\mathbf{H1}-\mathbf{H3}$ hold, then, according to
Proposition 3.1 in \cite{BDGM2014}, $\lambda (\theta )$ is the spectral
radius of $P_{\theta }$ and $P_{\theta }^{\ast }$ and there exist a unique
strictly positive function $r_{\theta }\in \mathcal{C}_{b}\left( \mathbb{X}%
\right) $ and a unique probability measure $l_{\theta }$ such that
\begin{equation*}
\int_{\mathbb{X}}r_{\theta }(\mathbf{x})dl_{\theta }(\mathbf{x})=1
\end{equation*}%
and%
\begin{equation}
l_{\theta }P_{\theta }=\lambda \left( \theta \right) l_{\theta },\quad
P_{\theta }r_{\theta }=\lambda \left( \theta \right) r_{\theta }.
\label{invariance}
\end{equation}%
Here $l_{\theta }P_{\theta }$ is a unique probability measure satisfying the
equality
\begin{equation*}
\int_{\mathbb{X}}g\left( \mathbf{x}\right) (l_{\theta }P_{\theta })\left( d%
\mathbf{x}\right) =\int_{\mathbb{X}}(P_{\theta }g\left( \mathbf{x}\right)
)l_{\theta }\left( d\mathbf{x}\right)
\end{equation*}%
for each $g\in \mathcal{C}_{b}\left( \mathbb{X}\right) $ (see \cite{CM2016},
page 2070) and $P_{\theta }r_{\theta }=\lambda \left( \theta \right)
r_{\theta }$ means that
\begin{equation}
\mathbb{E}\left[ \left\vert \mathbf{Mx}\right\vert ^{\theta }r_{\theta
}\left( \mathbf{M\circ x}\right) \right] =\lambda \left( \theta \right)
r_{\theta }(\mathbf{x})  \label{EigenLambda}
\end{equation}%
for all $\mathbf{x\in }\mathbb{X}$.

Similarly, there exists a pair $\left( r_{\theta }^{\ast },l_{\theta }^{\ast
}\right) $ possessing the same properties relative to $P_{\theta }^{\ast }$
as $\left( r_{\theta },l_{\theta }\right) $ relative to $P_{\theta }$.
Moreover,%
\begin{equation*}
r_{\theta }(\mathbf{x})=c\int_{\mathbb{X}}\left( \mathbf{x,y}\right)
^{\theta }l_{\theta }^{\ast }\left( d\mathbf{y}\right)
\end{equation*}%
where%
\begin{equation*}
c^{-1}=\int_{\mathbb{X\times X}}\left( \mathbf{x,y}\right) ^{\theta
}l_{\theta }^{\ast }\left( d\mathbf{x}\right) l_{\theta }\left( d\mathbf{y}%
\right) .
\end{equation*}

Following \cite{CM2016}, we introduce the functions
\begin{equation}
p_{n}^{(\theta )}\left( \mathbf{x},\mathbf{h}\right) :=\frac{\left\vert
\mathbf{hx}\right\vert ^{\theta }}{\lambda ^{n}\left( \theta \right) }\frac{%
r_{\theta }\left( \mathbf{h\circ x}\right) }{r_{\theta }\left( \mathbf{x}%
\right) },\quad \mathbf{x}\in \mathbb{X},n\geq 1.  \label{p-def}
\end{equation}

Using (\ref{EigenLambda}) it is not difficult to check that, for every $%
n\geq 1$, every $\mathbf{x}\in \mathbb{X}$ and every matrix $\mathbf{h\in }%
\mathbb{S}^{+}$
\begin{equation}
\mathbb{E}p_{n}^{(\theta )}\left( \mathbf{x},\mathbf{Mh}\right)
=p_{n-1}^{(\theta )}\left( \mathbf{x},\mathbf{h}\right) .
\label{consistency}
\end{equation}%
In particular,
\begin{equation}
\mathbb{E}p_{n}^{(\theta )}\left( \mathbf{x},\mathbf{L}_{n,1}\right) =1,
\label{total_mass}
\end{equation}%
where $\mathbf{L}_{n,k}:=\mathbf{M}_{n}\cdot \cdot \cdot \mathbf{M}_{k}$ if $%
1\leq k\leq n$ and $\mathbf{L}_{n,n+1}$ is the unit $p\times p$ matrix.

Thus, given $\mathbf{x}\in \mathbb{X}$ the function $p_{n}^{(\theta )}\left(
\mathbf{x},\mathbf{h}\right) $ may be viewed as the density of a probability
distribution on $\mathbb{S}^{+}$.

For each $n\geq 1$ denote by $\mathcal{F}_{n}$ the $\sigma $-algebra
generated by random elements $\mathbf{Z}_{1},\mathbf{Z}_{2},\ldots ,\mathbf{Z%
}_{n}$ and $\mathbf{F}_{1},\mathbf{F}_{2},\ldots ,\mathbf{F}_{n}$. Let $%
\mathbb{I}_{\mathcal{A}}$ be the indicator of the event $\mathcal{A}$.\ \ It
follows from \eqref{total_mass} that
\begin{equation}
\mathbb{P}_{n}^{(\theta )}(\mathcal{A})=\mathbb{P}_{n,\mathbf{x}}^{(\theta
)}(\mathcal{A}):=\mathbb{E}\left[ p_{n}^{(\theta )}\left( \mathbf{x},\mathbf{%
L}_{n,1}\right) \mathbb{I}_{\mathcal{A}}\right] ,\ \mathcal{A}\in \mathcal{F}%
_{n},  \label{measure0}
\end{equation}%
is a probability measure on $\mathcal{F}_{n}$. Furthermore, %
\eqref{consistency} implies that the sequence of measures $\left\{ \mathbb{P}%
_{n}^{(\theta )},n\geq 1\right\} $ is consistent and can be extended to a
probability measure $\mathbb{P}^{(\theta )}$ on our original probability
space $(\Omega ,\mathcal{F})$.

\section{Asymptotic of the survival probability\label{Sec2}}

We prove in this section that $\mathbb{P}\left( \mathcal{D}_{n}\left(
\mathbf{z}\right) \right) \sim K(\mathbf{z})\lambda ^{n}(1)$ as $%
n\rightarrow \infty $ for some constant $K(\mathbf{z})$. The precise
expression for $K(\mathbf{z})$ will be found in Section \ref{SecPrecise}.

For every environment $\mathcal{E}=$ $\left\{ \mathbf{F}_{n},n\geq 1\right\}
$ and $0\leq k<n$ we define iterations
\begin{eqnarray*}
&&\mathbf{F}_{k,n}(\mathbf{s}):=\mathbf{F}_{k+1}(\mathbf{F}_{k+2}(\ldots
\mathbf{F}_{n}(\mathbf{s})...))=\left( F_{k,n}^{1}(\mathbf{s}%
),F_{k,n}^{2}\left( \mathbf{s}\right) ,...,F_{k,n}^{p}(\mathbf{s})\right) ,
\\
&&\mathbf{F}_{n,k}(\mathbf{s}):=\mathbf{F}_{n}(\mathbf{F}_{n-1}(\ldots
\mathbf{F}_{k+1}(\mathbf{s})...))=\left( F_{n,k}^{1}(\mathbf{s}%
),F_{n,k}^{2}\left( \mathbf{s}\right) ,...,F_{n,k}^{p}(\mathbf{s})\right)
\end{eqnarray*}%
and set
\begin{equation*}
\mathbf{F}_{n,n}(\mathbf{s}):=\mathbf{s}.
\end{equation*}%
By definition of the process $\mathbf{Z}_{n}$
\begin{equation*}
\mathbb{E}\left[ \mathbf{s}^{\mathbf{Z}_{n}}\Big|\mathbf{Z}_{0}=\mathbf{e}%
_{i},\mathbf{F}_{1},\mathbf{F}_{2},\ldots ,\mathbf{F}_{n}\right]
=F_{0,n}^{i}(\mathbf{s}).
\end{equation*}%
Therefore,
\begin{equation}
\mathbb{P}\left( \mathbf{Z}_{n}\neq \mathbf{0}|\mathbf{Z}_{0}=\mathbf{e}_{i},%
\mathbf{F}_{1},\mathbf{F}_{2},\ldots ,\mathbf{F}_{n}\right) \mathbb{=}1%
\mathbb{-}F_{0,n}^{i}(\mathbf{0}).  \label{start.point.1}
\end{equation}%
This equality and definition (\ref{DefB}) imply
\begin{equation}
\mathbb{P}\left( \mathcal{D}_{n}\left( \mathbf{z}\right) \right) =\mathbb{E}%
\left[ 1-\prod_{i=1}^{p}\left( F_{0,n}^{i}(\mathbf{0})\right) ^{z_{i}}\right]
.  \label{ReprB}
\end{equation}

Let $\overline{\theta }:=(\theta _{1},...,\theta _{p})$ be a vector with
nonnegative real components $\theta _{i}\in \left\{ 0\right\} \cup \lbrack
1,\infty ),i=1,...,p,$ such that $\theta =\theta _{1}+...+\theta _{p}\geq 1$%
. Along with (\ref{ReprB}) we consider the quantity%
\begin{equation}
Y(n,\overline{\theta }):=\mathbb{E}\left[ \prod_{i=1}^{p}\left(
1-F_{0,n}^{i}(\mathbf{0})\right) ^{\theta _{i}}\right] =\mathbb{E}\left[
\prod_{i=1}^{p}\left( 1-F_{n,0}^{i}(\mathbf{0})\right) ^{\theta _{i}}\right]
.  \label{ReprAY}
\end{equation}

This function has the following probabilistic meaning. Let $\mathbf{Z}%
_{n}^{ij}$ is the vector of the number of particles in generation $n$ that
are offsprings of the $j$th particle of type $i$ existing in the population
at time $0$. Note that if $\overline{\theta }\mathbf{=z\in }\mathbb{N}%
_{+}^{p}$ and
\begin{equation}
\mathcal{A}_{n}\left( \mathbf{z}\right) :=\bigcap {}_{i=1}^{p}\bigcap
{}_{j=1}^{z_{i}}\left\{ \mathbf{Z}_{n}^{ij}\neq \mathbf{0}\big|\mathbf{Z}%
_{0}^{ij}=\mathbf{e}_{i}\right\}  \label{DefA}
\end{equation}%
is the event that each initial particle has a nonempty number of descendants
at moment $n$, then
\begin{equation}
Y(n,\mathbf{z})=\mathbb{P}\left( \mathcal{A}_{n}\left( \mathbf{z}\right)
\right) =\mathbb{E}\left[ \prod_{i=1}^{p}\left( 1-F_{0,n}^{i}(\mathbf{0}%
)\right) ^{z_{i}}\right] .  \label{ReprA}
\end{equation}

Since
\begin{equation}
1-F_{n,0}^{i}(\mathbf{s})=(\mathbf{e}_{i},\mathbf{1}-\mathbf{F}_{n,0}(%
\mathbf{s}))\leq (\mathbf{e}_{i},\mathbf{M}_{n}\left( \mathbf{1}-\mathbf{F}%
_{n-1,0}(\mathbf{s})\right) ),  \label{DefFi}
\end{equation}%
it follows that%
\begin{equation*}
1-F_{n,0}^{i}(\mathbf{0})\leq \mathbf{e}_{i}\mathbf{L}_{n,1}\mathbf{1}\leq |%
\mathbf{L}_{n,1}\mathbf{1}|.
\end{equation*}%
Using this representation we get by (\ref{ReprAY})
\begin{equation}
Y(n,\overline{\theta })\leq \mathbb{E}\left[ \prod_{i=1}^{p}|\mathbf{L}_{n,1}%
\mathbf{1}|^{\theta _{i}}\right] =\mathbb{E}\left[ |\mathbf{L}_{n,1}\mathbf{1%
}|^{\theta }\right] ,  \label{repr1}
\end{equation}%
where $\theta =\theta _{1}+...+\theta _{p}$.

\begin{lemma}
\label{L_arbit}Let $\overline{\theta }:=(\theta _{1},...,\theta _{p})$ be a
vector with nonnegative components $\theta _{i}\in \left\{ 0\right\} \cup
\lbrack 1,\infty ),i=1,...,p,$ such that $\theta =\left\vert \overline{%
\theta }\right\vert =\theta _{1}+...+\theta _{p}\geq 1$. Assume that
Conditions $\mathbf{H1}-\mathbf{H3}$ are valid and $\theta \in \Theta $.
Then there exists a constant $C=C(\theta )$ such that
\begin{equation*}
Y(n,\overline{\theta })\leq C\lambda ^{n}(\theta )
\end{equation*}%
for all $n\geq 1$.
\end{lemma}

\textbf{Proof}. We make the change of measure at the right-hand side of (\ref%
{repr1}) by means of the density $p_{n}^{(\theta )}\left( \mathbf{x},\mathbf{%
h}\right) $ in (\ref{p-def}) with $\theta =|\overline{\theta }|$ and $%
\mathbf{x}=p^{-1}\mathbf{1}$. As a result we get
\begin{eqnarray}
Y(n,\overline{\theta }) &\leq &\lambda ^{n}(\theta )r_{\theta }(\mathbf{x})%
\mathbb{E}\left[ p_{n}^{(\theta )}(\mathbf{x},\mathbf{L}_{n,1})\frac{1}{%
r_{\theta }(\mathbf{L}_{n,1}\circ \mathbf{x})}\frac{|\mathbf{L}_{n,1}\mathbf{%
1}|^{\theta }}{|\mathbf{L}_{n,1}\mathbf{x}|^{\theta }}\right]  \notag \\
&=&p^{\theta }\lambda ^{n}(\theta )r_{\theta }(p^{-1}\mathbf{1})\mathbb{E}%
^{(\theta )}\left[ \frac{1}{r_{\theta }(\mathbf{L}_{n,1}\circ \left( p^{-1}%
\mathbf{1}\right) )}\right] .  \label{ChangeOfMeasure}
\end{eqnarray}

Since $r_{\theta }(\mathbf{x})$ is a continuous positive function on the
compact $\mathbb{X}$, there are\ constants $c_{1}$ and $c_{2}$ such that
\begin{equation}
0<c_{1}\leq r_{\theta }(\mathbf{x})\leq c_{2}<\infty  \label{Rseparated}
\end{equation}%
for all $\mathbf{x}\in \mathbb{X}$. Hence it follows that
\begin{equation*}
Y(n,\overline{\theta })\leq p^{\theta }\frac{c_{2}}{c_{1}}\lambda
^{n}(\theta )=:C\lambda ^{n}(\theta ).
\end{equation*}

Lemma \ref{L_arbit} is proved.

Having this result in hands we check the validity of the following statement.

\begin{lemma}
\label{L_survivAll} Let Conditions $\mathbf{H1}-\mathbf{H5}$ be valid, $1\in
\mathrm{int}\Theta $ and $\Lambda ^{\prime }(1)<0$.Then for any $\mathbf{z}%
\in \mathbb{N}_{+}^{p}$, as $n\rightarrow \infty $%
\begin{equation}
\mathbb{P}\left( \mathcal{D}_{n}\left( \mathbf{z}\right) \right) \sim
\sum_{i=1}^{p}z_{i}\mathbb{P}\left( \mathcal{D}_{n}\left( \mathbf{e}%
_{i}\right) \right) .  \label{AsymEquiv}
\end{equation}
\end{lemma}

\textbf{Remark 3}. One may think that equivalence (\ref{AsymEquiv}) is
evident, since $\mathbb{P}\left( \mathcal{D}_{n}\left( \mathbf{e}_{i}\right)
\right) \rightarrow 0$ as $n\rightarrow \infty $. Indeed, this is always the
case for the ordinary multitype Galton-Watson processes. However, such
equivalence is, in general, not true for BPRE's. For instance, Bansaye \cite%
{Ban2009} has shown, analyzing a single type subcritical BPRE $\left\{
Z_{n},n\geq 0\right\} $, that if the process is weakly subcritical (i.e.,
satisfies\ the conditions $\Lambda ^{\prime }(0)<0,$ $\Lambda ^{\prime
}(1)>0 $ in our notation) then there exists $\beta \in \left( 0,1\right) $
such that the quantity
\begin{equation*}
\alpha _{z}:=\lim_{n\rightarrow \infty }\frac{\mathbb{P}\left(
Z_{n}>0|Z_{0}=z\right) }{\mathbb{P}\left( Z_{n}>0|Z_{0}=1\right) }
\end{equation*}%
is of order $z^{\beta }\log z$ as $z\rightarrow \infty $.

\textbf{Proof of Lemma \ref{L_survivAll}}. We need the following
inequalities
\begin{equation}
0\leq \sum_{k=1}^{p}z_{k}\left( b_{k}-a_{k}\right) -\left(
\prod_{k=1}^{p}b_{k}^{z_{k}}-\prod_{k=1}^{p}a_{k}^{z_{k}}\right) \leq
\sum_{k,l=1}^{p}z_{k}z_{l}\left( b_{k}-a_{k}\right) \left( b_{l}-a_{l}\right)
\label{DoublIneq}
\end{equation}%
valid for $z_{k}\in \mathbb{N}_{0}$ and real numbers $0\leq a_{k}\leq
b_{k}\leq 1,k=1,2,...,p$. Hence, using~(\ref{ReprB}) we obtain that, for any
$\theta _{0}\in \left( 1,2\right) $%
\begin{eqnarray}
0\leq\sum_{i=1}^{p}z_{i}\mathbb{P}\left( \mathcal{D}_{n}\left( \mathbf{e}%
_{i}\right) \right) -\mathbb{P}\left( \mathcal{D}_{n}\left( \mathbf{z}%
\right) \right) \leq \sum_{i,j=1}^{p}z_{i}z_{j}\mathbb{E}\left[ \left(
1-F_{0,n}^{i}(\mathbf{0})\right) \left( 1-F_{0,n}^{j}(\mathbf{0})\right) %
\right]  \notag \\
\leq2\left\vert \mathbf{z}\right\vert \sum_{i,j=1}^{p}z_{i}\mathbb{E}\left[
\left( 1-F_{0,n}^{i}(\mathbf{0})\right) ^{2}\right] \leq 2\left\vert \mathbf{%
z}\right\vert \sum_{i,j=1}^{p}z_{i}\mathbb{E}\left[ \left( 1-F_{0,n}^{i}(%
\mathbf{0})\right) ^{\theta _{0}}\right] .  \label{Negl2}
\end{eqnarray}

Since $\Lambda ^{\prime }\left( 1\right) <0$ and $1\in \mathrm{int}\Theta $,
there exists $\theta _{0}>1$ such that $\theta _{0}\in \mathrm{int}\Theta $
and $\Lambda ^{\prime }\left( \theta _{0}\right) <0$. By Lemma \ref{L_arbit}
\begin{equation}
\mathbb{E}\left[ \left( 1-F_{0,n}^{i}(\mathbf{0})\right) ^{\theta _{0}}%
\right] \leq C\lambda ^{n}(\theta _{0})  \label{Negl3}
\end{equation}%
as $n\rightarrow \infty $. Since $\lambda (\theta )$ is a convex function
for $\theta \in \Theta $, it follows that $\lambda (1)>\lambda (\theta
_{0})>0 $. Combining (\ref{Negl2}) and (\ref{SinglStrong0}) gives, as $%
n\rightarrow \infty $%
\begin{equation}
\mathbb{P}\left( \mathcal{D}_{n}\left( \mathbf{z}\right) \right) \sim
\sum_{i=1}^{p}z_{i}\mathbb{P}\left( \mathcal{D}_{n}\left( \mathbf{e}%
_{i}\right) \right) \sim \sum_{i=1}^{p}z_{i}K(\mathbf{e}_{i})\lambda ^{n}(1)
\label{SSingle}
\end{equation}%
with still not specified explicitly positive constants $K(\mathbf{e}_{i})$.

Lemma \ref{L_survivAll} is proved.

\section{Yaglom limit theorem\label{Sec3}}

Now we extend statement (\ref{EqYz0}) to the case of processes initiated by
several particles.

\begin{theorem}
\label{T_condYaglom} Assume that Conditions $\mathbf{H1}-\mathbf{H5}$ are
valid, $1\in \mathrm{int}\Theta $ and $\Lambda ^{\prime }(1)<0$. Then
\begin{equation}
\lim_{n\rightarrow \infty }\mathbb{E}\left[ \mathbf{s}^{\mathbf{Z}_{n}}|%
\mathbf{Z}_{n}\neq \mathbf{0};\mathbf{Z}_{0}=\mathbf{z}\right] =\sum_{%
\mathbf{j}\in \mathbb{N}_{+}^{p}}t_{\mathbf{j}}(\mathbf{z})\mathbf{s}^{%
\mathbf{j}}=:T\left( \mathbf{z},\mathbf{s}\right)
\end{equation}%
exists for all $\mathbf{z\in }\mathbb{N}_{+}^{p}$, specifies a proper
distribution on $\mathbb{N}_{+}^{p}$ and solves the equation
\begin{equation}
\mathbb{E}\left[ T\left( \mathbf{z},\mathbf{F}(\mathbf{s})\right) \right]
=\lambda (1)T\left( \mathbf{z},\mathbf{s}\right) +1-\lambda (1).
\label{Uniqueness0}
\end{equation}
\end{theorem}

\textbf{Proof}. We write%
\begin{eqnarray*}
&&\mathbb{E}\left[ \mathbf{s}^{\mathbf{Z}_{n}}|\mathbf{Z}_{n}\neq \mathbf{0};%
\mathbf{Z}_{0}=\mathbf{z}\right] =\frac{\mathbb{E}\left[ \mathbf{s}^{\mathbf{%
Z}_{n}};\mathbf{Z}_{n}\neq \mathbf{0}|\mathbf{Z}_{0}=\mathbf{z}\right] }{%
\mathbb{P}\left( \mathbf{Z}_{n}\neq \mathbf{0}|\mathbf{Z}_{0}=\mathbf{z}%
\right) } \\
&&\qquad \qquad =1-\frac{\mathbb{E}\left[ 1-\mathbf{s}^{\mathbf{Z}_{n}}|%
\mathbf{Z}_{0}=\mathbf{z}\right] }{\mathbb{P}\left( \mathbf{Z}_{n}\neq
\mathbf{0}|\mathbf{Z}_{0}=\mathbf{z}\right) }=1-\frac{\mathbb{E}\left[
1-\prod_{i=1}^{p}\left( F_{0,n}^{i}(\mathbf{s})\right) ^{z_{i}}\right] }{%
\mathbb{P}\left( \mathcal{D}_{n}\left( \mathbf{z}\right) \right) }.
\end{eqnarray*}%
By (\ref{DoublIneq}), (\ref{Negl2}), and (\ref{Negl3}) we see that%
\begin{eqnarray}
0 &\leq &\sum_{i=1}^{p}z_{i}\mathbb{E}\left[ 1-F_{0,n}^{i}(\mathbf{s})\right]
-\mathbb{E}\left[ 1-\prod_{i=1}^{p}\left( F_{0,n}^{i}(\mathbf{s})\right)
^{z_{i}}\right]  \notag \\
&\leq &\sum_{i,j=1}^{p}z_{i}z_{j}\mathbb{E}\left[ \left( 1-F_{0,n}^{i}(%
\mathbf{s})\right) \left( 1-F_{0,n}^{j}(\mathbf{s})\right) \right]  \notag \\
&\leq &2\left\vert \mathbf{z}\right\vert \sum_{i,j=1}^{p}z_{i}\mathbb{E}%
\left[ \left( 1-F_{0,n}^{i}(\mathbf{0})\right) ^{\theta _{0}}\right]
=o\left( \lambda ^{n}(1)\right)  \label{BuBu1}
\end{eqnarray}%
as $n\rightarrow \infty $. As shown in Theorem 1 of \cite{VW2018}
\begin{equation}
\lim_{n\rightarrow \infty }\frac{\mathbb{E}\left[ 1-F_{0,n}^{i}(\mathbf{s})%
\right] }{\mathbb{P}\left( \mathcal{D}_{n}\left( \mathbf{e}_{i}\right)
\right) }=:\phi _{i}(\mathbf{s})  \label{Wah}
\end{equation}%
exists for any $i\in \left\{ 1,...,p\right\} ,$ where $\phi _{i}(\mathbf{s}),%
\mathbf{s}\mathbf{\in \lbrack }0,1\mathbf{]}^{p},$ is a nondegenerate
multivariate function with $\phi _{i}(\mathbf{1})=0$.

Using (\ref{SSingle}) and setting%
\begin{equation*}
\mathbf{K}^{\ast }=\left( K\left( \mathbf{e}_{1}\right) ,...,K\left( \mathbf{%
e}_{p}\right) \right) ,
\end{equation*}%
where $K\left( \mathbf{e}_{i}\right) ,i=1,...,p,$ are the same as in (\ref%
{SinglStrong0}), we see that \
\begin{equation*}
\lim_{n\rightarrow \infty }\frac{\mathbb{P}\left( \mathcal{D}_{n}\left(
\mathbf{e}_{i}\right) \right) }{\mathbb{P}\left( \mathcal{D}_{n}\left(
\mathbf{z}\right) \right) }=\frac{K\left( \mathbf{e}_{i}\right) }{(\mathbf{z}%
,\mathbf{K}^{\ast })}.
\end{equation*}%
Thus, as $n\rightarrow \infty $%
\begin{eqnarray*}
\frac{1}{\mathbb{P}\left( \mathcal{D}_{n}\left( \mathbf{z}\right) \right) }%
\sum_{i=1}^{p}z_{i}\mathbb{E}\left[ 1-F_{0,n}^{i}(\mathbf{s})\right]
&=&\sum_{i=1}^{p}z_{i}\frac{\mathbb{P}\left( \mathcal{D}_{n}\left( \mathbf{e}%
_{i}\right) \right) }{\mathbb{P}\left( \mathcal{D}_{n}\left( \mathbf{z}%
\right) \right) }\frac{\mathbb{E}\left[ 1-F_{0,n}^{i}(\mathbf{s})\right] }{%
\mathbb{P}\left( \mathcal{D}_{n}\left( \mathbf{e}_{i}\right) \right) } \\
&\sim &\frac{1}{(\mathbf{z},\mathbf{K}^{\ast })}\sum_{i=1}^{p}z_{i}K\left(
\mathbf{e}_{i}\right) \phi _{i}(\mathbf{s}).
\end{eqnarray*}%
Setting $\Phi _{i}(\mathbf{s})=1-\phi _{i}(\mathbf{s})$ we obtain
\begin{eqnarray}
T\left( \mathbf{z},\mathbf{s}\right) := &&\lim_{n\rightarrow \infty }\mathbb{%
E}\left[ \mathbf{s}^{\mathbf{Z}_{n}}|\mathbf{Z}_{n}\neq \mathbf{0};\mathbf{Z}%
_{0}=\mathbf{z}\right]  \notag \\
&=&\frac{\sum_{i=1}^{p}z_{i}K\left( \mathbf{e}_{i}\right) \Phi _{i}(\mathbf{s%
})}{(\mathbf{z},\mathbf{K}^{\ast })}=\sum_{\mathbf{j}\in \mathbb{N}%
_{+}^{p}}t_{\mathbf{j}}(\mathbf{z})\mathbf{s}^{\mathbf{j}},  \label{DefYz1}
\end{eqnarray}%
where $T\left( \mathbf{z},\mathbf{s}\right) $ is a nondegenerate
multivariate function in $\mathbf{s}$ and $T\left( \mathbf{z},\mathbf{1}%
\right) =1$. Since the generating function $\mathbf{F}_{n+1}$ is independent
of the tuple $\mathbf{F}_{1},...,\mathbf{F}_{n}$ we get by (\ref{Wah}) and
the dominated convergence theorem that, as $n\rightarrow \infty $
\begin{eqnarray*}
&&\frac{\mathbb{E}\left[ 1-F_{0,n+1}^{i}(\mathbf{s})\right] }{\mathbb{P}%
\left( \mathcal{D}_{n+1}\left( \mathbf{z}\right) \right) }=\frac{\mathbb{P}%
\left( \mathcal{D}_{n}\left( \mathbf{e}_{i}\right) \right) }{\mathbb{P}%
\left( \mathcal{D}_{n+1}\left( \mathbf{z}\right) \right) }\frac{\mathbb{E}%
\left[ 1-F_{0,n}^{i}(\mathbf{F}_{n+1}(\mathbf{s}))\right] }{\mathbb{P}\left(
\mathcal{D}_{n}\left( \mathbf{e}_{i}\right) \right) } \\
&&\qquad =\frac{\mathbb{P}\left( \mathcal{D}_{n}\left( \mathbf{e}_{i}\right)
\right) }{\mathbb{P}\left( \mathcal{D}_{n+1}\left( \mathbf{z}\right) \right)
}\int_{\mathbf{y}\in \left[ 0,1\right] ^{p}}\frac{\mathbb{E}\left[
1-F_{0,n}^{i}(\mathbf{y})\right] }{\mathbb{P}\left( \mathcal{D}_{n}\left(
\mathbf{e}_{i}\right) \right) }\mathbb{P}\left( \mathbf{F}_{n+1}(\mathbf{s}%
)\in d\mathbf{y}\right) \\
&&\qquad \rightarrow \frac{1}{\lambda (1)}\frac{K(\mathbf{e}_{i})}{(\mathbf{z%
},\mathbf{K}^{\ast })}\mathbb{E}\left[ \phi _{i}(\mathbf{F}_{n+1}(\mathbf{s}%
))\right] =\frac{1}{\lambda (1)}\frac{K(\mathbf{e}_{i})}{(\mathbf{z},\mathbf{%
K}^{\ast })}\mathbb{E}\left[ \phi _{i}(\mathbf{F}(\mathbf{s}))\right] .
\end{eqnarray*}%
Therefore,%
\begin{eqnarray*}
&&T\left( \mathbf{z},\mathbf{s}\right) =\lim_{n\rightarrow \infty }\mathbb{E}%
\left[ \mathbf{s}^{\mathbf{Z}_{n+1}}|\mathbf{Z}_{n+1}\neq \mathbf{0};\mathbf{%
Z}_{0}=\mathbf{z}\right] \\
&=&1-\lim_{n\rightarrow \infty }\sum_{i=1}^{p}z_{i}\frac{\mathbb{E}\left[
1-F_{0,n+1}^{i}(\mathbf{s})\right] }{\mathbb{P}\left( \mathcal{D}%
_{n+1}\left( \mathbf{z}\right) \right) } \\
&=&1-\frac{1}{\lambda (1)}\sum_{i=1}^{p}z_{i}\frac{K(\mathbf{e}_{i})}{(%
\mathbf{z},\mathbf{K}^{\ast })}\mathbb{E}\left[ \phi _{i}(\mathbf{F}(\mathbf{%
s}))\right] =1-\frac{1}{\lambda (1)}\mathbb{E}\left[ 1-T\left( \mathbf{z},%
\mathbf{F}(\mathbf{s})\right) \right] .
\end{eqnarray*}%
Thus,%
\begin{equation}
\mathbb{E}\left[ T\left( \mathbf{z},\mathbf{F}(\mathbf{s})\right) \right]
=\lambda (1)T\left( \mathbf{z},\mathbf{s}\right) +1-\lambda (1).
\label{Uniqueness}
\end{equation}

Theorem \ref{T_condYaglom} is proved.

\section{Second change of measure\label{SecPrecise}}

To get more information about the limiting function $T\left( \mathbf{z},%
\mathbf{s}\right) $ we need one more measure ${\mathbb{P}^{\ast }}$ on the $%
\sigma $-field generated by $\mathbf{Z}_{1},\mathbf{Z}_{2},\ldots ;\mathbf{F}%
_{1},\mathbf{F}_{2},\ldots $. To this aim we fix $\mathbf{Z}_{0}=\mathbf{z}$
and, for every non-negative measurable functional $\varphi $ on $\mathbb{N}%
_{+}^{kp}\times \mathcal{P}_{p}^{k}(\mathbb{N}_{0}^{p}),\,k\geq 1,$ set
\begin{equation}
{\mathbb{E}^{\ast }}[\varphi (\mathbf{Z}_{1},\ldots ,\mathbf{Z}_{k};\mathbf{F%
}_{1},\ldots ,\mathbf{F}_{k})]\,:=\,\frac{\mathbb{E}[\varphi (\mathbf{Z}%
_{1},\ldots ,\mathbf{Z}_{k};\mathbf{F}_{1},\ldots ,\mathbf{F}_{k})(\mathbf{Z}%
_{k},\mathbf{U})]}{(\mathbf{z},\mathbf{U})\lambda }.  \label{size2b}
\end{equation}

\begin{lemma}
\label{L_consist} Relation (\ref{size2b}) defines a probability measure ${%
\mathbb{P}^{\ast }}$ on the $\sigma $--algebra generated by $\mathbf{Z}_{1},%
\mathbf{Z}_{2},\ldots ;\mathbf{F}_{1},\mathbf{F}_{2},\ldots $.
\end{lemma}

\textbf{Proof}. Set $\mathbf{R}_{1,k}=\mathbf{M}_{1}\mathbf{M}_{2}\cdot
\cdot \cdot \mathbf{M}_{k}$. Clearly,
\begin{eqnarray*}
\mathbb{E}\left[ \mathbf{Z}_{k}|\mathbf{Z}_{0}=\mathbf{z}\right] &=&\mathbb{E%
}\left[ \mathbb{E}[\mathbf{Z}_{k}\mid \mathbf{Z}_{1},\ldots ,\mathbf{Z}%
_{k-1};\mathcal{E}]\right] \\
&=&\mathbb{E}\left[ \mathbf{Z}_{k-1}\mathbf{M}_{k}|\mathbf{Z}_{0}=\mathbf{z}%
\right] =\mathbb{E}[\mathbf{zM}_{1}\cdot \cdot \cdot \mathbf{M}_{k}] \\
&=&\mathbf{z}\mathbb{E}[\mathbf{R}_{1,k}]=\mathbf{zm}^{k}.
\end{eqnarray*}%
Hence, using (\ref{MuEigen}) we see that
\begin{eqnarray*}
{\mathbb{E}^{\ast }}[1] &:&=\frac{\mathbb{E}[(\mathbf{Z}_{k},\mathbf{U})|%
\mathbf{Z}_{0}=\mathbf{z}]}{(\mathbf{z},\mathbf{U})\lambda ^{k}}=\frac{(%
\mathbb{E}[\mathbf{Z}_{k}|\mathbf{Z}_{0}=\mathbf{z}],\mathbf{U})}{(\mathbf{z}%
,\mathbf{U})\lambda ^{k}} \\
&=&\frac{(\mathbf{zm}^{k},\mathbf{U})}{(\mathbf{z},\mathbf{U})\lambda ^{k}}=%
\frac{(\mathbf{z},\mathbf{m}^{k}\mathbf{U})}{(\mathbf{z},\mathbf{U})\lambda
^{k}}=1.
\end{eqnarray*}%
Thus, $\mathbb{P}{^{\ast }}$ is, indeed, a probability measure for each $k$.
Furthermore, the following consistency condition holds: If $\mathbf{Z}_{0}=%
\mathbf{z}\in \mathbb{N}_{+}^{p}$ and functions $\varphi _{k}$ and $\varphi
_{k+1}$ satisfy
\begin{equation*}
\varphi _{k+1}(\mathbf{z}_{1},\ldots ,\mathbf{z}_{k+1};\mathbf{f}_{1},\ldots
,\mathbf{f}_{k+1})=\varphi _{k}(\mathbf{z}_{1},\ldots ,\mathbf{z}_{k};%
\mathbf{f}_{1},\ldots ,\mathbf{f}_{k})
\end{equation*}%
for all $\mathbf{z}_{i}\in \mathbb{N}_{0}^{p}$ and $\mathbf{f}_{i}\in
\mathcal{P}_{p}(\mathbb{N}_{0}^{p}),\,1\leq i\leq k+1,$ then
\begin{align*}
& \mathbb{E}^{\ast }[\varphi _{k+1}(\mathbf{Z}_{1},\ldots ,\mathbf{Z}_{k+1};%
\mathbf{F}_{1},\ldots ,\mathbf{F}_{k+1})] \\
& =\frac{\mathbb{E}[\varphi _{k+1}(\mathbf{Z}_{1},\ldots ,\mathbf{Z}_{k};%
\mathbf{F}_{1},\ldots ,\mathbf{F}_{k+1})(\mathbf{Z}_{k+1},\mathbf{U})]}{(%
\mathbf{z},\mathbf{U})\lambda ^{k+1}} \\
& =\frac{\mathbb{E}[\varphi _{k}(\mathbf{Z}_{1},\ldots ,\mathbf{Z}_{k};%
\mathbf{F}_{1},\ldots ,\mathbf{F}_{k})(\mathbf{Z}_{k+1},\mathbf{U})]}{(%
\mathbf{z},\mathbf{U})\lambda ^{k+1}} \\
& =\frac{\mathbb{E}\big[\varphi _{k}(\mathbf{Z}_{1},\ldots ,\mathbf{Z}_{k};%
\mathbf{F}_{1},\ldots ,\mathbf{F}_{k})\mathbb{E}[(\mathbf{Z}_{k+1},\mathbf{U}%
)\mid \mathbf{Z}_{1},\ldots ,\mathbf{Z}_{k};\mathcal{E}]\big]}{(\mathbf{z},%
\mathbf{U})\lambda ^{k+1}} \\
& =\frac{\mathbb{E}\big[\varphi _{k}(\mathbf{Z}_{1},\ldots ,\mathbf{Z}_{k};%
\mathbf{F}_{1},\ldots ,\mathbf{F}_{k})(\mathbf{Z}_{k}\mathbb{E}\mathbf{M}%
_{k+1},\mathbf{U})]\big]}{(\mathbf{z},\mathbf{U})\lambda ^{k+1}} \\
& =\frac{\mathbb{E}\big[\varphi _{k}(\mathbf{Z}_{1},\ldots ,\mathbf{Z}_{k};%
\mathbf{F}_{1},\ldots ,\mathbf{F}_{k})(\mathbf{Z}_{k},\mathbf{mU})]\big]}{(%
\mathbf{z},\mathbf{U})\lambda ^{k+1}} \\
& =\frac{\mathbb{E}[\varphi _{k}(\mathbf{Z}_{1},\ldots ,\mathbf{Z}_{k};%
\mathbf{F}_{1},\ldots ,\mathbf{F}_{k})(\mathbf{Z}_{k},\mathbf{U})]}{(\mathbf{%
z},\mathbf{U})\lambda ^{k}} \\
& =\mathbb{E}^{\ast }[\varphi _{k}(\mathbf{Z}_{1},\ldots ,\mathbf{Z}_{k};%
\mathbf{F}_{1},\ldots ,\mathbf{F}_{k})],
\end{align*}%
where for the last three equalities we have used the relation
\begin{equation*}
\mathbb{E}\left[ \mathbf{Z}_{k+1}\mid \mathbf{Z}_{k}=\mathbf{x},\mathcal{E}=(%
\mathbf{f}_{1},\mathbf{f}_{2},\ldots )\right] =\mathbf{xM}_{k+1},
\end{equation*}%
and the independency of $\mathbf{M}_{k+1}$ and the random vector $(\mathbf{Z}%
_{1},\ldots ,\mathbf{Z}_{k},\mathbf{F}_{1},\ldots ,\mathbf{F}_{k})$.

The lemma is proved.

\textbf{Remark 4}. Note that the measure $\mathbb{P}^{\ast }$ has
essentially different nature in comparison with the measure $\mathbb{P}%
^{(1)} $ defined by formula (\ref{measure0}). For instance, under $\mathbb{P}%
^{(1)}$ the MBPRE $\mathbf{Z}_{n}$ may take value $\mathbf{0}$ with
probability
\begin{equation*}
\mathbb{P}_{n}^{(1)}(\mathbf{Z}_{n}=\mathbf{0})=\mathbb{E}\left[
p_{n}^{(1)}\left( \mathbf{x},\mathbf{L}_{n,1}\right) \mathbb{I}_{\left\{
\mathbf{Z}_{n}=\mathbf{0}\right\} }\right] >0,
\end{equation*}%
while under $\mathbb{P}^{\ast }$ the probability of the event $\left\{
\mathbf{Z}_{n}=\mathbf{0}\right\} $ always is zero.

\textbf{Remark 5}. If $\varphi _{k}(\mathbf{z}_{1},\ldots ,\mathbf{z}_{k};%
\mathbf{f}_{1},\ldots ,\mathbf{f}_{k})=\varphi _{k}(\mathbf{f}_{1},\ldots ,%
\mathbf{f}_{k})$, i.e. depends only on generating functions, then
\begin{align}
\mathbb{E}^{\ast }[\varphi (\mathbf{F}_{1},\ldots ,\mathbf{F}_{k})]& =\frac{%
\mathbb{E}\big[\varphi (\mathbf{F}_{1},\ldots ,\mathbf{F}_{k})\mathbb{E}[(%
\mathbf{Z}_{k},\mathbf{U})\mid \mathcal{E}]\big]}{(\mathbf{z},\mathbf{U}%
)\lambda ^{k}}  \notag \\
& =\frac{\mathbb{E}\big[\varphi (\mathbf{F}_{1},\ldots ,\mathbf{F}%
_{k})\left( \mathbb{E[}\mathbf{Z}_{k-1}\mid \mathcal{E}]\mathbf{M}_{k},%
\mathbf{U}\right) \big]}{(\mathbf{z},\mathbf{U})\lambda ^{k}}  \notag \\
& =\frac{\mathbb{E}[\varphi (\mathbf{F}_{1},\ldots ,\mathbf{F}_{k})(\mathbf{%
zM}_{1}\cdot \cdot \cdot \mathbf{M}_{k},\mathbf{U})]}{(\mathbf{z},\mathbf{U}%
)\lambda ^{k}}  \notag \\
& =\frac{\mathbb{E}[\varphi (\mathbf{F}_{1},\ldots ,\mathbf{F}_{k})(\mathbf{%
zR}_{1,k},\mathbf{U})]}{(\mathbf{z},\mathbf{U})\lambda ^{k}}  \label{size2}
\end{align}%
for each $k\in \mathbb{N}$.

We recall that according to the definition of MBPRE's the transition
probabilities $P_{\mathbf{xy}}$ of the Markov chain $\{\mathbf{Z}_{n},n\geq
0\}$ are
\begin{equation}
P_{\mathbf{xy}}:=\mathbb{E}[(\mathbf{F}^{\mathbf{x}})[\mathbf{y}]],\quad
\mathbf{x},\mathbf{y}\in \mathbb{N}_{0}^{p},  \label{original}
\end{equation}%
where $(\mathbf{f}^{\mathbf{x}})[\mathbf{y}]$ is the weight assigned by $%
\mathbf{f}^{\mathbf{x}}$ to the point $\mathbf{y}\in \mathbb{N}_{0}^{p}$.
Using this definition we can identify the distribution of the process~$\{%
\mathbf{Z}_{n},n\geq 0\}$ under the measure $\mathbb{P}^{\ast }$ as the law
of a certain Markov chain. Note that (\ref{size2b}) implies
\begin{align}
\mathbb{P}^{\ast }(\mathbf{Z}_{1}=\mathbf{z}_{1},\ldots ,\mathbf{Z}_{k}=%
\mathbf{z}_{k}|\mathbf{Z}_{0}=\mathbf{z}_{0})& :=\frac{(\mathbf{z}_{k},%
\mathbf{U})\,\mathbb{P}(\mathbf{Z}_{1}=\mathbf{z}_{1},\ldots ,\mathbf{Z}_{k}=%
\mathbf{z}_{k}|\mathbf{Z}_{0}=\mathbf{z}_{0})}{(\mathbf{z}_{0},\mathbf{U}%
)\lambda ^{k}}  \notag \\
& =\frac{(\mathbf{z}_{1},\mathbf{U})\,\cdots (\mathbf{z}_{k},\mathbf{U})\,}{(%
\mathbf{z}_{0},\mathbf{U})\cdots (\mathbf{z}_{k-1},\mathbf{U})}\,\frac{P_{%
\mathbf{z}_{0}\mathbf{z}_{1}}\cdots P_{\mathbf{z}_{k-1}\mathbf{z}_{k}}}{%
\lambda ^{k}}\   \notag \\
& =\ \prod_{j=1}^{k}P_{\mathbf{z}_{j-1}\mathbf{z}_{j}}^{\ast }
\label{Schachtel11}
\end{align}%
for every $\mathbf{z}_{j}\in \mathbb{N}_{+}^{p},\,0\leq j\leq k,$ where
\begin{equation}
P_{\mathbf{xy}}^{\ast }\,:=\,\frac{(\mathbf{y},\mathbf{U})P_{\mathbf{xy}}}{(%
\mathbf{x},\mathbf{U})\lambda },\ \ \mathbf{x},\mathbf{y}\in \mathbb{N}%
_{+}^{p}.  \label{kanada}
\end{equation}

By linearity of expectation%
\begin{eqnarray*}
\sum_{\mathbf{y}\in \mathbb{N}_{+}^{p}}P_{\mathbf{xy}}^{\ast }\, &=&\frac{1}{%
(\mathbf{x},\mathbf{U})\lambda }\sum_{\mathbf{y}\in \mathbb{N}_{+}^{p}}(%
\mathbf{y},\mathbf{U})P_{\mathbf{xy}}=\frac{1}{(\mathbf{x},\mathbf{U}%
)\lambda }\mathbb{E}\left[ (\mathbf{Z}_{1},\mathbf{U})|\mathbf{Z}_{0}=%
\mathbf{x}\right] \\
&=&\frac{1}{(\mathbf{x},\mathbf{U})\lambda }(\mathbf{x}\mathbb{E}\left[
\mathbf{M}_{1}\right] ,\mathbf{U})=\frac{1}{(\mathbf{x},\mathbf{U})\lambda }(%
\mathbf{x,\mathbf{m}U})=1
\end{eqnarray*}%
for every $\mathbf{x}\in \mathbb{N}_{+}^{p}$ . Thus, $P^{\ast }$ indeed,
specifies a Markov chain on $\mathbb{N}_{+}^{p}$.

There is an evident connection between the $k$-step transition probabilities
of the $P$-chain and the $k$-step transition probabilities of the $P^{\ast }$%
-chain which we express by the following formula:
\begin{equation}
\mathbb{P}^{\ast }(\mathbf{Z}_{j+k}=\mathbf{y}\,|\,\mathbf{Z}_{j}=\mathbf{x}%
)\,=\,P_{\mathbf{xy}}^{\ast k}=\frac{(\mathbf{y},\mathbf{U})P_{\mathbf{xy}%
}^{k}}{(\mathbf{x},\mathbf{U})\lambda ^{k}}  \label{kanada2}
\end{equation}%
for every $\mathbf{x},\mathbf{y}\in \mathbb{N}_{+}^{p}$ and $j,k\geq 0$.
Here $P^{\ast k},P^{k}$ denote the $k$th power of the transition matrices $%
P^{\ast }$ and $P$.

\ For \ $\mathbf{x=}\left( x_{1},...,x_{p}\right) \in \mathbb{N}_{+}^{p}$ we
write the representation
\begin{eqnarray}
\mathbb{E}\left[ \left( \mathbf{F}(\mathbf{s})\right) ^{\mathbf{x}}\right]
&=&\mathbb{E}\left( \prod_{i=1}^{p}\left( F^{i}(\mathbf{s})\right)
^{x_{i}}\right)  \notag \\
&=&\sum_{\mathbf{y}\in \mathbb{N}_{0}^{p}}P_{\mathbf{xy}}\mathbf{s}^{\mathbf{%
y}} =\sum_{\mathbf{y}\in \mathbb{N}_{0}^{p}}\mathbb{P}\left( \mathbf{Z}_{1}=%
\mathbf{y}|\mathbf{Z}_{0}=\mathbf{x}\right) \mathbf{s}^{\mathbf{y}}.
\label{Expansion}
\end{eqnarray}%
Thus,
\begin{eqnarray*}
\mathbb{E}\left[ T\left( \mathbf{z},\mathbf{F}(\mathbf{s})\right) \right]
&=&\sum_{\mathbf{x}\in \mathbb{N}_{+}^{p}}t_{\mathbf{x}}(\mathbf{z})\mathbb{E%
}\left( \mathbf{F}(\mathbf{s})\right) ^{\mathbf{x}} \\
&=&\sum_{\mathbf{x}\in \mathbb{N}_{+}^{p}}t_{\mathbf{x}}(\mathbf{z})\sum_{%
\mathbf{y}\in \mathbb{N}_{0}^{p}}P_{\mathbf{xy}}\mathbf{s}^{\mathbf{y}%
}=\sum_{\mathbf{y}\in \mathbb{N}_{0}^{p}}\mathbf{s}^{\mathbf{y}}\sum_{%
\mathbf{x}\in \mathbb{N}_{+}^{p}}t_{\mathbf{x}}(\mathbf{z})P_{\mathbf{xy}}.
\end{eqnarray*}%
Using (\ref{Uniqueness0}) we get
\begin{equation*}
\sum_{\mathbf{y}\in \mathbb{N}_{0}^{p}}\mathbf{s}^{\mathbf{y}}\sum_{\mathbf{x%
}\in \mathbb{N}_{+}^{p}}t_{\mathbf{x}}(\mathbf{z})P_{\mathbf{xy}}=\lambda
\sum_{\mathbf{y}\in \mathbb{N}_{+}^{p}}t_{\mathbf{y}}(\mathbf{z})\mathbf{s}^{%
\mathbf{y}}+1-\lambda .
\end{equation*}%
Thus, if $\mathbf{y}\neq \mathbf{0}$ then%
\begin{equation}
\sum_{\mathbf{x}\in \mathbb{N}_{+}^{p}}t_{\mathbf{x}}(\mathbf{z})P_{\mathbf{%
xy}}=\lambda t_{\mathbf{y}}(\mathbf{z}),  \label{PreStationar}
\end{equation}%
and if $\mathbf{y}=\mathbf{0}$ then%
\begin{equation}
\sum_{\mathbf{x}\in \mathbb{N}_{+}^{p}}t_{\mathbf{x}}(\mathbf{z})P_{\mathbf{%
x0}}=1-\lambda .  \label{Remainder}
\end{equation}

Note that by Fatou's lemma
\begin{eqnarray}
W(\mathbf{z}):= &&\sum_{\mathbf{y}\in \mathbb{N}_{+}^{p}}(\mathbf{y},\mathbf{%
U})t_{\mathbf{y}}(\mathbf{z})=\sum_{\mathbf{y}\in \mathbb{N}_{+}^{p}}(%
\mathbf{y},\mathbf{U})\lim_{n\rightarrow \infty }\mathbb{P}\left( \mathbf{Z}%
_{n}=\mathbf{y}|\mathbf{Z}_{n}\neq \mathbf{0;Z}_{0}=\mathbf{z}\right)
\notag \\
&\leq &\liminf_{n\rightarrow \infty }\mathbb{E}\left[ (\mathbf{Z}_{n},%
\mathbf{U})|\mathbf{Z}_{n}\neq \mathbf{0;Z}_{0}=\mathbf{z}\right]   \notag \\
&=&\liminf_{n\rightarrow \infty }\frac{\mathbb{E}\left[ (\mathbf{Z}_{n},%
\mathbf{U})|\mathbf{Z}_{0}=\mathbf{z}\right] }{\mathbb{P}\left( \mathbf{Z}%
_{n}\neq \mathbf{0|Z}_{0}=\mathbf{z}\right) }  \notag \\
&=&\liminf_{n\rightarrow \infty }\frac{\left( \mathbb{E}\left[ \mathbf{Z}%
_{n}|\mathbf{Z}_{0}=\mathbf{z}\right] ,\mathbf{U}\right) }{\mathbb{P}\left(
\mathbf{Z}_{n}\neq \mathbf{0|Z}_{0}=\mathbf{z}\right) }=\lim_{n\rightarrow
\infty }\frac{\left( \mathbf{z}\mathbb{E}\left[ \mathbf{R}_{1,n}\right] ,%
\mathbf{U}\right) }{\mathbb{P}\left( \mathbf{Z}_{n}\neq \mathbf{0|Z}_{0}=%
\mathbf{z}\right) }  \notag \\
&=&\lim_{n\rightarrow \infty }\frac{\left( \mathbf{z,m}^{n}\mathbf{U}\right)
}{\mathbb{P}\left( \mathbf{Z}_{n}\neq \mathbf{0|Z}_{0}=\mathbf{z}\right) }%
=\lim_{n\rightarrow \infty }\frac{\lambda ^{n}\left( \mathbf{z,U}\right) }{(%
\mathbf{z},\mathbf{K}^{\ast })\lambda ^{n}}=\frac{\left( \mathbf{z,U}\right)
}{(\mathbf{z},\mathbf{K}^{\ast })}<\infty .  \label{FiniteMean}
\end{eqnarray}%
Let
\begin{equation}
t_{\mathbf{x}}^{\ast }(\mathbf{z})=\frac{\left( \mathbf{x,U}\right) t_{%
\mathbf{x}}(\mathbf{z})}{W(\mathbf{z})},\ \mathbf{x}\in N_{+}^{p}.
\label{DefStat}
\end{equation}%
Clearly, \ for each \ $\mathbf{z}\in \mathbb{N}_{+}^{p}$
\begin{equation}
\sum_{\mathbf{x}\in \mathbb{N}_{+}^{p}}t_{\mathbf{x}}^{\ast }(\mathbf{z})=1.
\label{Sum22}
\end{equation}%
Using (\ref{kanada}) we obtain for $\mathbf{y}\in \mathbb{N}_{+}^{p}$
\begin{eqnarray*}
\sum_{\mathbf{x}\in \mathbb{N}_{+}^{p}}t_{\mathbf{x}}^{\ast }(\mathbf{z})P_{%
\mathbf{xy}}^{\ast } &=&\sum_{\mathbf{x}\in \mathbb{N}_{+}^{p}}\frac{\left(
\mathbf{x},\mathbf{U}\right) t_{\mathbf{x}}(\mathbf{z})}{W(\mathbf{z})}%
\times \frac{\left( \mathbf{y},\mathbf{U}\right) P_{\mathbf{xy}}}{\left(
\mathbf{x},\mathbf{U}\right) \lambda } \\
&=&\frac{\left( \mathbf{y},\mathbf{U}\right) }{W(\mathbf{z})\lambda }\sum_{%
\mathbf{x}\in \mathbb{N}_{+}^{p}}t_{\mathbf{x}}(\mathbf{z})P_{\mathbf{xy}}=%
\frac{\left( \mathbf{y},\mathbf{U}\right) }{W(\mathbf{z})}t_{\mathbf{y}}(%
\mathbf{z})=t_{\mathbf{y}}^{\ast }(\mathbf{z}).
\end{eqnarray*}%
As a result,
\begin{equation}
\sum_{\mathbf{x}\in \mathbb{N}_{+}^{p}}t_{\mathbf{x}}^{\ast }(\mathbf{z})P_{%
\mathbf{xy}}^{\ast }=t_{\mathbf{y}}^{\ast }\left( \mathbf{z}\right) ,\ \
\mathbf{y}\in \mathbb{\ N}_{+}^{p}.  \label{invariance2}
\end{equation}

\begin{lemma}
\label{L_invariant}Let Conditions $\mathbf{H1}-\mathbf{H6}$ be valid, $1\in
\mathrm{int}\Theta $ and $\Lambda ^{\prime }(1)<0$. Then

(i) The probability measure $t^{\ast }$ from (\ref{DefStat}) is an invariant
distribution for the Markov chain generated by transition probabilities $P_{%
\mathbf{xy}}^{\ast }$.

(ii)The chain has a single recurrent class $R=\mathrm{supp}\,t^{\ast }$. The
class~$R$ is positive recurrent and aperiodic.
\end{lemma}

\textbf{Proof.} (i) By (\ref{Sum22}), the measure $t^{\ast }$ has total mass
1. The invariance follows from (\ref{invariance2}).

(ii) We first show that there are states which can be reached from any other
state of the chain in two steps. Recall that \ $m_{ij}>0$ for all $i,j\in
\left\{ 1,...,p\right\} $. Therefore, for any $\mathbf{x}\in \mathbb{N}%
_{+}^{p}$ and $i\in \left\{ 1....,p\right\} $%
\begin{equation*}
\sum_{\mathbf{w=}\left( w_{1},...,w_{p}\right) \in \mathbb{N}%
_{+}^{p}:w_{i}\geq 1}P_{\mathbf{xw}}>0.
\end{equation*}%
Further, in view of Condition\textbf{\ }$\mathbf{H6}$, there exist $i\in
\left\{ 1,...,p\right\} $ and $\mathbf{y=}\left( y_{1},...,y_{p}\right) \in
\mathbb{N}_{+}^{p}$ with $y_{i}\geq 1$ such that
\begin{equation*}
\mathbb{P}(F^{i}[\mathbf{y}]\neq \mathbf{0}\text{ and }F^{k}[\mathbf{0}%
]>0,k=1,...,p)>0,
\end{equation*}%
i.e., type $i$ individuals of the same generation of the original branching
process may have both $\mathbf{0}$ or $\mathbf{y\neq 0}$ children with a
positive probability. For such $\mathbf{y}$ and any \ $\mathbf{x}\in \mathbb{%
N}_{+}^{p}$
\begin{equation*}
P_{\mathbf{xy}}^{2}\geq \sum_{\mathbf{w}\in \mathbb{N}_{+}^{p}:w_{i}\geq
1}P_{\mathbf{xw}}P_{\mathbf{wy}}\geq \sum_{\mathbf{w}\in \mathbb{N}%
_{+}^{p}:w_{i}\geq 1}P_{\mathbf{xw}}\mathbb{E}\big[F^{i}[\mathbf{y}](\mathbf{%
F}[\mathbf{0}])^{\mathbf{w}-\mathbf{e}_{i}}\big]>0
\end{equation*}%
and, therefore,
\begin{equation}
P_{\mathbf{xy}}^{\ast 2}=\,\frac{(\mathbf{y},\mathbf{U})P_{\mathbf{xy}}^{2}}{%
(\mathbf{x},\mathbf{U})\lambda ^{2}}>0.  \label{positive}
\end{equation}%
Besides, for the \ $\mathbf{y}$
\begin{equation}
P_{\mathbf{yy}}^{\ast }=\lambda ^{-1}P_{\mathbf{yy}}\geq \lambda ^{-1}%
\mathbb{E}\big[F^{i}[\mathbf{y}](\mathbf{F}[\mathbf{0}])^{\mathbf{y}-\mathbf{%
e}_{i}}\big]>0.  \label{positive2}
\end{equation}

The second assertion of the lemma now follows from standard results from
Markov chain theory: Since any invariant probability distribution is
supported by positive recurrent states (see, e.g. the criterion in
Section~XV.7 of~\cite{Fel1968}), part (i) of the proposition shows that the
chain has at least one such class. In view of (\ref{positive}) there can be
at most one recurrent class. Clearly, this class~$R$, say, contains all $%
\mathbf{y}$ which satisfy~(\ref{positive2}). Since $P_{\mathbf{yy}}^{\ast
}>0 $ for such $\mathbf{y}$, the class is aperiodic. The fact that $R=%
\mbox{supp}\,t^{\ast }$ again follows from part (i), because the equilibrium
weight $t_{\mathbf{y}}^{\ast }(\mathbf{z})$ is the reciprocal of the
expected return time to~$\mathbf{y}$ (see, e.g., Theorem~1 in Section~XV.7
of~\cite{Fel1968}) and, therefore, is unique for all $\mathbf{z}\in \mathbb{N%
}_{+}^{p}$.

The lemma is proved.

\begin{lemma}
\label{L_UNiq} If the conditions of Theorem \ref{T_condYaglom} are valid
then the function $T\left( \mathbf{z},\mathbf{s}\right) $ in (\ref%
{Uniqueness0}) is one and the same for all $\mathbf{z\in }\mathbb{N}_{+}^{p}$%
.
\end{lemma}

\textbf{Proof}. Since $P^{\ast }$ has a unique recurrent class and does not
depend on a particular $\mathbf{z}$ it follows that the quantities $t_{%
\mathbf{y}}^{\ast }(\mathbf{z}),\mathbf{y}\in \mathbb{N}_{+}^{p},$ also do
not depend on $\mathbf{z}$. We write them as $t_{\mathbf{y}}^{\ast }$.
Recalling (\ref{DefStat}) we see that, for each $\mathbf{y}$ the ratio
\begin{equation}
\frac{t_{\mathbf{y}}^{\ast }}{\left( \mathbf{y},\mathbf{U}\right) }=\frac{t_{%
\mathbf{y}}(\mathbf{z})}{W(\mathbf{z})}  \label{Ratio}
\end{equation}%
is one and the same for all $\mathbf{z}\in \mathbb{N}_{+}^{p}$. This, in
turn, means that the left-hand side of the equality
\begin{equation*}
\sum_{\mathbf{y}\in \mathbb{N}_{+}^{p}}\frac{t_{\mathbf{y}}(\mathbf{z})}{W(%
\mathbf{z})}P_{\mathbf{y0}}=\frac{1-\lambda }{W(\mathbf{z})}
\end{equation*}%
is one and the same for all $\mathbf{z}\in \mathbb{N}_{+}^{p}$. Therefore, $%
W(\mathbf{z})=W$ is also independent of $\mathbf{z}$. In view of (\ref{Ratio}%
) the last implies independence of $t_{\mathbf{y}}(\mathbf{z})=t_{\mathbf{y}%
} $ of $\mathbf{z}\in \mathbb{N}_{+}^{p}$. \ Therefore, $T\left( \mathbf{z},%
\mathbf{s}\right) $ is one and the same for all $\mathbf{z\in }\mathbb{N}%
_{+}^{p}\mathbf{,}$ and we write%
\begin{equation*}
T\left( \mathbf{s}\right) =\sum_{\mathbf{y}\in \mathbb{N}_{+}^{p}}t_{\mathbf{%
y}}\mathbf{s}^{\mathbf{y}}=\mathbb{E}\left[ \mathbf{s}^{\mathbf{T}_{\infty }}%
\right] ,\quad T\left( \mathbf{1}\right) =1,
\end{equation*}%
where, in view of (\ref{FiniteMean})%
\begin{equation*}
W=\sum_{\mathbf{y}\in \mathbb{N}_{+}^{p}}\left( \mathbf{y},\mathbf{U}\right)
t_{\mathbf{y}}=\mathbb{E}\left[ \left( \mathbf{T}_{\infty },\mathbf{U}%
\right) \right] <\infty ,
\end{equation*}%
and%
\begin{equation*}
t_{\mathbf{y}}^{\ast }=\frac{\left( \mathbf{y},\mathbf{U}\right) t_{\mathbf{y%
}}}{W},\mathbf{y}\in N_{+}^{p}.
\end{equation*}%
This completes the proof of Lemma \ref{L_UNiq}.

Now we would like to find an explicit form for the vector $\mathbf{K}^{\ast
} $ in (\ref{SinglStrong}).

\begin{lemma}
\label{L_constant}Under the conditions of Theorem \ref{T_condSeveral}
\begin{equation*}
\mathbf{K}^{\ast }\mathbf{=}\frac{\mathbf{U}}{\mathbb{E}\left[ \left(
\mathbf{T}_{\infty },\mathbf{U}\right) \right] }.
\end{equation*}
\end{lemma}

\textbf{Proof}. It follows from (\ref{EqYz}) and (\ref{SinglStrong}) that
for any $\mathbf{z\in }\mathbb{N}_{+}^{p}$
\begin{equation*}
\lim_{n\rightarrow \infty }\frac{\mathbb{P}\left( \mathbf{Z}_{n}=\mathbf{y}%
\,|\mathbf{Z}_{0}=\mathbf{z}\right) }{\mathbb{P}\left( \mathbf{Z}_{n}\neq
\mathbf{0}\,|\mathbf{Z}_{0}=\mathbf{z}\right) }=\lim_{n\rightarrow \infty }%
\frac{P_{\mathbf{zy}}^{n}}{\mathbb{(}\mathbf{z,K}^{\ast }\mathbb{)}\lambda
^{n}}=t_{\mathbf{y}}.
\end{equation*}%
Given $t_{\mathbf{y}}>0$ we get
\begin{equation*}
t_{\mathbf{y}}^{\ast }=\frac{(\mathbf{y,U})}{(\mathbf{z},\mathbf{U})}t_{%
\mathbf{y}}=\lim_{n\rightarrow \infty }\frac{(\mathbf{y,U})P_{\mathbf{zy}%
}^{n}}{(\mathbf{z},\mathbf{U})\lambda ^{n}}=\frac{(\mathbf{y,U})\mathbb{(}%
\mathbf{z,K}^{\ast }\mathbb{)}}{(\mathbf{z},\mathbf{U})}t_{\mathbf{y}}=W%
\frac{\mathbb{(}\mathbf{z,K}^{\ast }\mathbb{)}}{(\mathbf{z},\mathbf{U})}t_{%
\mathbf{y}}^{\ast }.
\end{equation*}%
Thus,
\begin{equation*}
\frac{\mathbb{(}\mathbf{z,K}^{\ast }\mathbb{)}}{(\mathbf{z},\mathbf{U})}=%
\frac{1}{W}
\end{equation*}%
and, therefore, the left-hand side is one and the same for all $\mathbf{z\in
}\mathbb{N}_{+}^{p}$. Selecting sequentially $\mathbf{z}=\mathbf{e}%
_{i},i=1,...,p,$ gives $\mathbf{K}^{\ast }=W^{-1}\mathbf{U=K}$, where $%
\mathbf{K}$ is the same vector as in (\ref{DefK}).

Lemma \ref{L_constant} is proved.

\textbf{Proof of Theorem \ref{T_condSeveral}.} Combining Theorem \ref%
{T_condYaglom} and Lemmas \ref{L_UNiq} and \ref{L_constant} gives the proof
of Theorem \ref{T_condSeveral}.

\section{Q-processes\label{SecQproc}}

\textbf{Proof of Corollary \ref{C_Qprocess2}}. Relation (\ref{PrelMult})
easily follows from Theorem \ref{T_condSeveral}.

Thus, we prove (\ref{LimitMult}). For $0=n_{0}<n_{1}<....<n_{r}=n$ we write
\begin{eqnarray}
&&\mathbb{P}\left( \mathbf{Z}_{n_{1}}=\mathbf{j}_{1},\mathbf{Z}_{n_{2}}=%
\mathbf{j}_{2},...,\mathbf{Z}_{n_{r}}=\mathbf{j}_{r}|\mathbf{Z}_{n_{r}}\neq
\mathbf{0;Z}_{0}=\mathbf{j}_{0}\right)  \notag \\
&&\quad =\frac{\lambda ^{n_{r}}}{\mathbb{P}\left( \mathbf{Z}_{n_{r}}\neq
\mathbf{0|Z}_{0}=\mathbf{j}_{0}\right) }\prod_{i=1}^{r}\frac{\mathbb{P}%
\left( \mathbf{Z}_{n_{i}}=\mathbf{j}_{i}|\mathbf{Z}_{n_{i}-1}=\mathbf{j}%
_{i-1}\right) }{\lambda ^{n_{i}-n_{i-1}}}.  \label{NN}
\end{eqnarray}%
We know by Theorem \ref{T_condSeveral} that
\begin{align}
&\lim_{n\rightarrow \infty }\mathbb{P}\left( \mathbf{Z}_{n}=\mathbf{y|Z}_{0}=%
\mathbf{x}\right)  \notag \\
&\qquad=\lim_{n\rightarrow \infty }\frac{\mathbb{P}\left( \mathbf{Z}_{n}\neq
0|\mathbf{Z}_{0}=\mathbf{x}\right) }{\lambda ^{n}}\mathbb{P}\left( \mathbf{Z}%
_{n}=\mathbf{y|Z}_{n}\neq 0,\mathbf{Z}_{0}=\mathbf{x}\right)  \notag \\
&\qquad=\left( \mathbf{y,K}\right) \mathbb{P}\left( \mathbf{T}_{\infty }=%
\mathbf{y}\right) .  \label{RR}
\end{align}%
Using (\ref{RR}) and (\ref{SinglStrong}) we deduce by (\ref{NN}) that%
\begin{align*}
\lim_{\min_{i}(n_{i}-n_{i-1)}\rightarrow \infty }\mathbb{P}\left( \mathbf{Z}%
_{n_{1}}=\mathbf{j}_{1},\mathbf{Z}_{n_{2}}=\mathbf{j}_{2},...,\mathbf{Z}_{n}=%
\mathbf{j}_{r}|\mathbf{Z}_{n_{r}}\neq \mathbf{0;Z}_{0}=\mathbf{j}_{0}\right)&
\\
=\frac{1}{\left( \mathbf{j}_{r}\mathbf{,K}\right) }\left(
\prod_{i=1}^{r}\left( \mathbf{j}_{i}\mathbf{,K}\right) \mathbb{P}\left(
\mathbf{T}_{\infty }=\mathbf{j}_{i}\right) \right)\qquad & \\
=\left( \prod_{i=1}^{r-1}\left( \mathbf{j}_{i}\mathbf{,K}\right) \mathbb{P}%
\left( \mathbf{T}_{\infty }=\mathbf{j}_{i}\right) \right) \mathbb{P}\left(
\mathbf{T}_{\infty }=\mathbf{j}_{r}\right)&
\end{align*}%
as desired.

\textbf{Proof of Corollary \ref{C_Qprocess}}. By definition
\begin{eqnarray*}
&&\mathbb{P}\left( \mathbf{\hat{Y}}_{1}=\mathbf{j}_{1},\mathbf{\hat{Y}}_{2}=%
\mathbf{j}_{2},...,\mathbf{\hat{Y}}_{n}=\mathbf{j}_{n}|\mathbf{\hat{Y}}_{0}=%
\mathbf{y}\right) \\
&&\qquad=\lim_{m\rightarrow \infty }\mathbb{P}\left( \mathbf{Z}_{1}=\mathbf{j%
}_{1},\mathbf{Z}_{2}=\mathbf{j}_{2},...,\mathbf{Z}_{n}=\mathbf{j}_{n}|%
\mathbf{Z}_{n+m}\neq \mathbf{0;Z}_{0}=\mathbf{y}\right) .
\end{eqnarray*}%
In view of (\ref{SinglStrong}) and Theorem \ref{T_condSeveral}
\begin{eqnarray*}
&&\mathbb{P}\left( \mathbf{\hat{Y}}_{1}=\mathbf{j}_{1},\mathbf{\hat{Y}}_{2}=%
\mathbf{j}_{2},...,\mathbf{\hat{Y}}_{n}=\mathbf{j}_{n}|\mathbf{\hat{Y}}_{0}=%
\mathbf{y}\right) \\
&&=\mathbb{P}\left( \mathbf{Z}_{1}=\mathbf{j}_{1},\mathbf{Z}_{2}=\mathbf{j}%
_{2},...,\mathbf{Z}_{n}=\mathbf{j}_{n}\mathbf{|Z}_{0}=\mathbf{y}\right)
\lim_{m\rightarrow \infty }\frac{\mathbb{P}\left( \mathbf{Z}_{m}\neq \mathbf{%
0|Z}_{0}=\mathbf{j}_{n}\right) }{\mathbb{P}\left( \mathbf{Z}_{n+m}\neq
\mathbf{0|Z}_{0}=\mathbf{y}\right) } \\
&&\qquad\qquad\qquad\quad=\frac{1}{\lambda ^{n}}\frac{\left( \mathbf{j}_{n}%
\mathbf{,K}\right) }{\left( \mathbf{y,K}\right) }\mathbb{P}\left( \mathbf{Z}%
_{1}=\mathbf{j}_{1},\mathbf{Z}_{2}=\mathbf{j}_{2},...,\mathbf{Z}_{n}=\mathbf{%
j}_{n}\mathbf{|Z}_{0}=\mathbf{y}\right) .
\end{eqnarray*}%
Hence, recalling (\ref{RR}) and (\ref{SinglStrong}) we get
\begin{align*}
\lim_{n\rightarrow \infty }\mathbb{P}\left( \mathbf{\hat{Y}}_{n}=\mathbf{j}|%
\mathbf{\hat{Y}}_{0}=\mathbf{y}\right) =\lim_{n\rightarrow \infty }\frac{1}{%
\lambda ^{n}}\frac{\left( \mathbf{j,K}\right) }{\left( \mathbf{y,K}\right) }%
\mathbb{P}\left( \mathbf{Z}_{n}=\mathbf{j,Z}_{n}\neq \mathbf{0|Z}_{0}=%
\mathbf{y}\right)& \\
=\left( \mathbf{j,K}\right) \lim_{n\rightarrow \infty }\frac{\mathbb{P}%
\left( \mathbf{Z}_{n}\neq \mathbf{0|Z}_{0}=\mathbf{y}\right) }{\left(
\mathbf{y,K}\right) \lambda ^{n}(1)}\mathbb{P}\left( \mathbf{Z}_{n}=\mathbf{%
j|Z}_{n}\neq \mathbf{0,Z}_{0}=\mathbf{y}\right)& \\
=\left( \mathbf{j,K}\right) \mathbb{P}\left( \mathbf{T}_{\infty }=\mathbf{j}%
\right) =\frac{\left( \mathbf{j,U}\right) }{\mathbb{E}\left[ \left( \mathbf{T%
}_{\infty },\mathbf{U}\right) \right] }\mathbb{P}\left( \mathbf{T}_{\infty }=%
\mathbf{j}\right).&
\end{align*}%
Corollary \ref{C_Qprocess2} is proved.

\end{document}